\documentclass[11pt]{amsart}
\usepackage{amsmath,amsthm, amscd, amssymb, amsfonts, mathrsfs}
\usepackage[all]{xy}
\usepackage{color}
\usepackage{enumitem}
\usepackage{hhline}
\usepackage[active]{srcltx}%{vista inversa}

\newcommand{\PSL}{\mathbf{PSL}}
\newcommand{\Gb}{\boldsymbol{G}}
\newcommand{\PSp}{\mathbf{PSp}}
\newcommand{\Pom}{\mathbf{P\Omega}}

\newcommand{\PSU}{\mathbf{PSU}}

\newcommand{\ba}{\mathbf{a}}
\newcommand{\bb}{\mathbf{b}}
\newcommand{\be}{\mathbf{e}}

\newcommand{\bu}{\mathbf{u}}
\newcommand{\bv}{\mathbf{v}}
\newcommand{\bw}{\mathbf{w}}

\newcommand{\x}{\mathtt{x}}
\newcommand{\y}{\mathtt{y}}

\newcommand{\St}{\mathtt{S}}
\newcommand{\Rt}{\mathtt{R}}

\newcommand\sn{\mathbb S_n}

\newcommand{\fd}{finite-dimensional}

\newcommand{\GCF}{\operatorname{gcf}}

\newcommand{\Ind}{\operatorname{Ind}}

\newcommand{\Irr}{\operatorname{Irr}}

\newcommand{\ord}{\operatorname{ord}}

\newcommand{\supp}{\operatorname{supp}}

\newcommand{\Inn}{\operatorname{Inn}}

\newcommand{\Tb}{\mathbb T}

\newcommand{\kc}{\mathbb F_q}

\newcommand{\kdos}{\mathbb F_2}

\newcommand\toba{\mathfrak B }
\newcommand\bpi{\boldsymbol{\pi}}
\newcommand{\gr}{\operatorname{gr}}
\newcommand{\trid}{\triangleright}

\newcommand{\Fc}{{\mathcal F}}

\newcommand{\W}{{\mathcal W}}
\newcommand{\Zc}{{\mathcal Z}}

\newcommand{\irrfact}{{\mathfrak{I}}}

\newcommand{\kk}{\Bbbk}%{}

\newcommand{\ku}{\mathbb C}
\newcommand{\Gc}{{\mathcal G}}
\newcommand{\K}{{\mathcal K}}

\newcommand{\Z}{{\mathbb Z}}
\newcommand{\N}{{\mathbb N}}
\newcommand{\I}{{\mathbb I}}
\newcommand{\J}{{\mathbb J}}

\newcommand{\G}{{\mathbb G}}

\newcommand{\T}{{\mathbb{T}}}
\newcommand{\U}{\mathbb{U}}

\newcommand{\F}{{\mathbb F}}
\newcommand{\C}{{\mathcal C}}
\newcommand{\D}{{\mathcal D}}

\newcommand{\GL}{\mathbf{GL}}
\newcommand{\PGL}{\mathbf{PGL}}
\newcommand{\SL}{\mathbf{SL}}

\newcommand{\uno}{{\bf 1}}

\newcommand{\Oc}{{\mathcal O}}
\newcommand{\oc}{{\mathcal O}}

\newcommand{\ydg}{{}^{\ku G}_{\ku G}\mathcal{YD}}

\newcommand\Gsc{\G_{\operatorname{sc}}}

\numberwithin{equation}{section}
\theoremstyle{plain}

\newtheorem{lema}{Lemma}[section]
\newtheorem{theorem}[lema]{Theorem}

\newtheorem{prop}[lema]{Proposition}

\newtheorem{question}{Question}
\newtheorem{question-app}{Question}

\theoremstyle{definition}
\newtheorem{definition}[lema]{Definition}
\newtheorem{exa}[lema]{Example}

\theoremstyle{remark}
\newtheorem{obs}[lema]{Remark}

\newtheorem{step}{Case}
\newtheorem{stepnss}{Case}

\newcommand\id{\operatorname{id}}

\newcommand\Id{\operatorname{id}}

\newcommand\st{\mathbb S_3}
\newcommand\sk{\mathbb S_4}

\newcommand\dn{\mathbb D_n}
\newcommand\dl{\mathbb D_d}

\newcommand\ac{\mathbb A_4}
\newcommand\aco{\mathbb A_5}

\newcommand\A{\mathbb A}
\newcommand\ao{\mathbb A_8}

\newcommand\Sim{\mathbb S}

\newcommand\s{\mathbb S}

\def\pf{\begin{proof}}
\def\epf{\end{proof}}

\theoremstyle{remark}
\newtheorem{punto}{}

\begin{document}

\renewcommand{\baselinestretch}{1.2}

\thispagestyle{empty}
%\vspace*{2in}
\title[Nichols algebras over non-semisimple classes in $\PSL_n(q)$]
{Finite-dimensional pointed Hopf algebras\newline over finite simple groups of Lie type \newline I. Non-semisimple classes in $\PSL_n(q)$}

\author[N. Andruskiewitsch, G. Carnovale, G. A. Garc\'ia]
{Nicol\'as Andruskiewitsch, Giovanna Carnovale and \\Gast\'on Andr\'es Garc\'ia}

\thanks{
This work was partially supported by
ANPCyT-Foncyt, CONICET, Secyt (UNC), the bilateral agreement between the Universities of C\'ordoba and Padova, the Visiting Scientist Program of the University
of Padova, the Grant CPDA 1258 18/12 of the University of Padova and the Visiting Professor Program of GNSAGA.}

\address{\noindent N. A. : Facultad de Matem\'atica, Astronom\'{\i}a y F\'{\i}sica,
Universidad Nacional de C\'ordoba. CIEM -- CONICET. %\newline \noindent
Medina Allende s/n (5000) Ciudad Universitaria, C\'ordoba,
Argentina}
\address{\noindent G. C. :
Dipartimento di Matematica,
Universit\`a degli Studi di Padova,
via Trieste 63, 35121 Padova, Italia}
\address{\noindent G. A. G. : Departamento de Matem\'atica, Facultad de Ciencias Exactas,
Universidad Nacional de La Plata. CONICET. Casilla de Correo 172, (1900)
La Plata, Argentina.}
\email{andrus@famaf.unc.edu.ar}
\email{carnoval@math.unipd.it}
\email{ggarcia@mate.unlp.edu.ar}

%\subjclass[2010]{16T05}
%\date{\today}

\begin{abstract}
We show that  Nichols algebras of most simple Yetter-Drin\-feld modules over the projective special linear group over a finite field,
corresponding to non-semisimple orbits, have infinite dimension. We spell out a new criterium to show that a rack collapses.
\end{abstract}

\maketitle

\begin{quote}{\textit{Ph'nglui mglw'nafh Cthulhu
R'lyeh wgah'nagl fhtagn}}\end{quote}

%\setcounter{tocdepth}{1}

%\tableofcontents

\section{Introduction}

This is the first article of a series  intended to determine the finite-dimensi\-onal
pointed Hopf algebras with group of group-likes isomorphic to a finite simple group of Lie type.
We now give an Introduction to the whole series.

\subsection{}
The general question we are dealing with is the classification
of \fd{} complex pointed Hopf algebras $H$ whose group of
group-like elements is a finite simple group. We say
that a finite group $G$ \emph{collapses} when every
finite-dimensional
pointed Hopf algebra $H$, with $G(H) \simeq G$ is
isomorphic to $\ku G$ \cite{AFGV-ampa}.
Here are the antecedents of that question.
\begin{itemize}
  \item If $G\simeq \Z/p$ is simple abelian, then the
  classification is known: for $p =2$ by \cite{N},
  see also \cite{CD}; for $p>7$, by
  \cite[Remark 1.10 (v)]{AS-ann-ens}; for $p = 5, 7$,
  combining  \cite[Theorem 1.3]{AS-adv} and \cite{AS-ann}.

\medbreak
  \item If $G\simeq \A_m$, $m\ge 5$ is alternating,
  then $G$ collapses \cite{AFGV-ampa}.

\medbreak
  \item If $G$ is a sporadic simple group, then $G$ collapses,
  except for the groups $G = Fi_{22}$, $B$, $M$. For these
groups, all irreducible Yetter--Drinfeld modules $M(\oc,\rho)$ have infinite
dimensional Nichols algebra, except for
a short list appearing in \cite[Table 1]{AFGV-espo} and improved in
\cite[Appendix]{FV}, of examples not known to be
finite-dimensional.

\medbreak
  \item If $G = \SL_2(q)$, $\GL_2(q)$, $\PGL_2(q)$ or
  $\PSL_2(q)$, all irreducible Yetter--Drinfeld modules
  $M(\oc,\rho)$ have infinite dimensional Nichols algebra,
  except for a list of examples not known to be
finite-dimensional given in \cite{FGV1, FGV2}. Particularly,
$\PSL_2(q)$ collapses for $q > 2$ even.
\end{itemize}

\subsection{}\label{subsec:intro-finite-Lie}
In this series of papers  we consider \fd{} pointed Hopf algebras with
finite simple group of Lie type.
We recall the description of such groups. Let $p$ be a prime
number, $m\in \N$, $q =p^m$ and $\kc$ the field with $q$ elements.

\medbreak
$\diamond$ Let $\G$ be a semisimple algebraic group
defined over  $\kc$.  A
Steinberg endomorphism $F: \G\to\G$ is an abstract group
automorphism having a power equal to a Frobenius map
\cite[21.3]{MT}.
The subgroup $\G^F$ of fixed points by $F$ is called a
\emph{finite group of Lie type} \cite[21.6]{MT}.

\medbreak
$\diamond$ Assume that $\G$ is a simple simply connected algebraic group.
Then $\G/Z(\G)$ is a simple abstract group \cite[12.5]{MT}
but $\G^F$ is not simple in general.
In fact $\Gb := \G^F/Z(\G^F)$ is a simple finite group
except for 8 examples that appear in low rank and with
$q =2$ or 3 (Tits Theorem, \cite[24.17]{MT}).
These $\Gb$ are called  \emph{finite simple groups of Lie type}
although they are \emph{not} finite groups of Lie type
in the sense above, in general.

\medbreak
There are three possible classes of Steinberg
endomorphisms of simple algebraic groups  \cite[22.5]{MT}
and accordingly we consider three families of
finite simple groups of Lie type:

\medbreak
\emph{Chevalley groups}. In the terminology of \cite{MT},
 these correspond to $\kc$-split Steinberg endomorphisms.
 That is, there exists an
$F$-stable torus $\Tb$ such that $F(t)=t^q$ for all $t\in \Tb$.
Then  $F$ is called a Frobenius map and $\G^F = \G(\kc)$
is the finite group of $\kc$-points.
Explicitly, these are the groups:

\medbreak $\PSL_n(q)$, $n\ge 2$ (except $\PSL_2(2)\simeq \st$
and $\PSL_2(3)\simeq \ac$ that are not simple);
$\PSp_{2n}(q)$, $n \ge 2$; $\Pom_{2n+1}(q)$, $n \ge 3$, $q$ odd;
$\Pom^+_{2n}(q)$, $n \ge 4$;
$G_2(q)$, $q\ge 3$; $F_4(q)$;
$E_6(q)$; $E_7(q)$; $E_8(q)$.

\medbreak
Also,
$\PSL_{2}(4)\simeq \PSL_{2}(5) \simeq \mathbb{A}_{5}$;
$\PSL_{2}(9) \simeq \mathbb{A}_{6}$;
 $\PSL_{4}(2) \simeq \mathbb{A}_{8}$, cf. \cite{W}.

\medbreak
\emph{Steinberg groups}.  These correspond to {\it twisted}
Steinberg endomorphisms.  A Steinberg endomorphism is twisted if it is not split and it is the
product of an $\F_q$-split endomorphism with an algebraic automorphism of $\G$ \cite{MT};
we may assume that $F$ is the product of a Frobenius
endomorphism with an automorphism of $\G$
induced by a non-trivial Dynkin diagram automorphism.
Explicitly, these are the groups:

\medbreak
$\PSU_n(q)$, $n \ge 3$ (except $\PSU_3(2)$);
$\Pom_{2n}^-(q)$, $n \ge 4$; ${}^3D_4(q)$, ${}^2E_6(q)$.

\medbreak
\emph{Suzuki-Ree groups}. Related to {\it very twisted}
Steinberg endomorphisms \cite{MT}.
Explicitly, these are the groups:

\medbreak
${}^2B_2(2^{2h+1})$, $h \ge 1$; ${}^2G_2(3^{2h+1})$,
$h \ge 1$; ${}^2F_4(2^{2h+1})$, $h \ge 1$.

\subsection{} The base field is $\ku$.
 Let $G$ be a finite group and let $H$ be a pointed Hopf algebra
 with $G(H) \simeq G$. For details on the following
 exposition-- not needed henceforth and included
 only for completeness, see \cite{AS-cambr, AG-adv}.

 Let $ 0 = H_{-1} \subset H_0 = \ku G(H) \subset H_1 \subset \dots$
 be the coradical filtration of $H$ and
 $\gr H = \oplus_{n\in \N_0} H_n/H_{n-1} \simeq R \# \ku G(H)$
 be the associated graded Hopf algebra.
 Here $R = \oplus_{n\in \N_0} R^n$ is a graded Hopf algebra in
 the braided tensor category $\ydg$
 of Yetter-Drinfeld modules over $\ku G$. Also, the subalgebra
 of $R$ generated by $V := R^1$
 is isomorphic to the Nichols algebra $\toba(V)$ of $V$. Hence
 $\dim H < \infty \iff \dim R < \infty \implies \dim \toba(V) <\infty$.
 Thus we need to address the  question:
 \emph{Determine all $V\in \ydg$ with $\dim \toba(V) < \infty$}.

\medbreak
In particular, the following are equivalent \cite[Lemma 1.4]{AFGV-ampa}:
\begin{itemize}
  \item $G$ collapses.

  \item For every $V\in \ydg$, $\dim \toba(V) = \infty$.

  \item For every \emph{irreducible} $V\in \ydg$, $\dim \toba(V) = \infty$.
\end{itemize}

Now all irreducible Yetter-Drinfeld modules over $\ku G$ are of
the form $M(\oc, \rho) = \Ind_{C_{G}(g)}^G \rho$,
where $\oc$ is a conjugacy class of $G$
and $\rho\in \Irr{C_{G}(g)}$ for $g\in \oc$ fixed. Set
$\toba(\oc, \rho) := \toba(M(\oc, \rho))$. Then the initial
question about the classification of \fd{} pointed Hopf
algebras with finite simple group of Lie type $G$ relies on the
consideration of the following:

\begin{question}\label{que:nichols-irrepyd-collapse} For such $G$,
determine all pairs $(\oc, \rho)$  with $\dim \toba(\oc, \rho) < \infty$.
\end{question}

\subsection{}\label{subsec:intro-collapse}
A crucial observation is that the algebra $\toba(\oc, \rho)$ does
not depend on the Yetter-Drinfeld module structure of $M(\oc, \rho)$
but only on the underlying braided vector space $(\ku \oc, c^{\rho})$.
In other words, the algebra $\toba(\oc, \rho)$ depends only on
the rack $\oc$ and the non-principal 2-cocycle arising from $\rho$,
see Section \ref{sec:racks} for definitions, or \cite{AG-adv} for more details.
In fact, to solve Question \ref{que:nichols-irrepyd-collapse} for
every finite group $G$ is tantamount to solve

\begin{question}\label{que:racks-cocycle-inf-dim}
\cite[Question 2]{AFGV-ampa}
Determine all pairs $(X, q)$, where $X$ is a finite rack and  $q$
is a non-principal 2-cocycle,
such that $\dim \toba(X, c^q) < \infty$.
\end{question}

The meaning of the next definition relies on the existence of
some criteria for a rack to collapse,
cf. \S \ref{subsec:intro-collapse-criteria}, \S \ref{subsec:typeD}.

\begin{definition}\label{def:rack-collapses} \cite[2.2]{AFGV-ampa}
A rack $X$ \emph{collapses} when $\dim \toba(X, q) = \infty$
for every finite faithful 2-cocycle $q$.
\end{definition}

Therefore, we tackle the initial question about the classification
of \fd{} pointed Hopf algebras with finite simple group of Lie type $G$
(rephrased as Question \ref{que:nichols-irrepyd-collapse})
in the following way:

$\bullet$ Determine all conjugacy classes in $G$ that collapse.

\smallbreak
$\bullet$ If $\oc$ is a conjugacy class in $G$ that does not collapse, then
  for any $\rho$ as above, compute the restriction $c^{\rho}_{X}$
  of the braiding $c^{\rho}$ to a suitable
  abelian subrack $X$ of $\oc$.
  If the Nichols algebra $\toba(\ku X, c^{\rho}_{X})$
  has infinite dimension (and this is checked by inspection of
  the list in \cite{H-all}),
  then so has $\toba(\ku \oc, c^{\rho})$.

\subsection{}\label{subsec:intro-collapse-criteria}
In principle, to solve Question \ref{que:racks-cocycle-inf-dim}
one would need first to compute all
possible non-principal 2-cocycles for a fixed rack $X$, before
starting to deal with the corresponding Nichols algebras.
A remarkable fact is the existence of criteria that dispense of this
computation. The first such criterium
is about racks of type D \cite{AFGV-ampa}, see \S \ref{subsec:typeD}:
\emph{If $X$ is a finite rack of type D, then $X$  collapses}.
In \S \ref{subsec:typeD} we introduce the notion
rack of type F, and prove an
analogous criterium. To distinguish the setting where neither of these criteria apply, 
we shall say that  a rack is \emph{cthulhu}\footnote{See \texttt{http://en.wikipedia.org/wiki/Cthulhu} for spelling and pronunciation.}
when it is neither of type D nor of type F. 
Also a rack is \emph{sober} if  every subrack is either abelian or indecomposable; this is  stronger than being cthulhu.
See \S \ref{subsec:cthulhu} for examples.

\subsection{} We need the description the conjugacy classes
in finite simple groups of Lie type.
Let $\G$ be a simple algebraic  group, $\Gsc$ its simply connected cover with $\bpi: \Gsc \to \G$ the natural projection, $F$ a Steinberg
endomorphism of $\Gsc$, cf. \S \ref{subsec:intro-finite-Lie},  $\Gb = \Gsc^F/Z(\Gsc^F)$,
 $\pi: \Gsc^F \to \Gb$ the natural projection. Often $F$ descends to $\G$, and then
there is a projection $\pi: [\G^F, \G^F] \to  [\G^F, \G^F]/ \bpi(Z(\Gsc^F)) \simeq \Gb$.
Every $x\in \Gsc$ has a Chevalley-Jordan decomposition
$x = x_sx_u= x_ux_s$, with $x_{s}$ semisimple and $x_{u}$ unipotent.
This decomposition boils down to $\G$ and to the finite groups $\Gsc^F$
and $\Gb$, where it agrees with the decomposition
in the $p$-part, namely $x_u$, and the $p$-regular part, namely $x_s$.
We state a well-known fact referred to as the \emph{isogeny argument}.
Let $\Gc$ be a semisimple algebraic, resp. finite, group and
$\Gc_u$ the set of unipotent, resp. $p$-elements, in $\Gc$.

\begin{lema}\label{lem:isogeny-argument}
Let $\Zc$ be a central (algebraic) subgroup of $\Gc$
consisting of semisimple, resp. $p$-regular elements.
Then the quotient map $\pi: \Gc \to \Gc/\Zc$
induces a rack isomorphism  $\pi: \Gc_u \to (\Gc/\Zc)_u$
and a bijection between the set of $\Gc$-conjugacy classes
in $\Gc_u$ and that of $\Gc/\Zc$-conjugacy classes in $(\Gc/\Zc)_u$.
\end{lema}

If $\Gc$ is semisimple algebraic, then $\Zc$ is finite because it
consists of semisimple elements. Hence $\Gc/\Zc$ is again a
semisimple algebraic group.

\pf Clearly $\pi(\Gc_u) \subset (\Gc/\Zc)_u$. Let $g\in \Gc$
with $\pi(g) \in (\Gc/\Zc)_u$ and let $g = g_sg_u$ be its Chevalley-Jordan
decomposition (resp., the decomposition in the $p$-regular
and the $p$-part). Then $\pi(g) = \pi(g_s) \pi(g_u)$, hence $\pi(g) = \pi(g_u)$
by uniqueness of the decomposition.
Thus $\pi: \Gc_u \to (\Gc/\Zc)_u$ is surjective.
Let now $g, h\in \Gc_u$ with $\pi(g) =\pi(h)$. Then $g =h z = zh$
for some $z\in \Zc$; but this turns out to be the decomposition of
$g$, hence $g = h$ and $\pi: \Gc_u \to (\Gc/\Zc)_u$ is a rack isomorphism.
Finally, let again $g, h\in \Gc_u$. If $\oc_g^{\Gc} = \oc_h^{\Gc}$,
then clearly $\oc_{\pi(g)}^{\Gc/\Zc} = \oc_{\pi(h)}^{\Gc/\Zc}$.
Conversely, if $\oc_{\pi(g)}^{\Gc/\Zc} = \oc_{\pi(h)}^{\Gc/\Zc}$,
then there exists $u\in \G$ and $z\in \Zc$ such that $ugu^{-1} = hz = zh$;
this is the decomposition of $ugu^{-1}\in \Gc_u$, hence $ugu^{-1} = h$
and $\oc_g^{\Gc} = \oc_h^{\Gc}$.
\epf

\medbreak
Let $x\in \Gb$; pick $\x\in \Gsc^F$ such that $\pi(\x) = x$.
If $\x = \x_s\x_u$ is
its Chevalley-Jordan decomposition, then
$x_s = \pi(\x_s)$, $x_u = \pi(\x_u)$ is the Chevalley-Jordan
decomposition of $x$, with $x_{s}$ semisimple and $x_{u}$ unipotent.
Now $\x_u$ belongs to $\K := C_{\Gsc^F}(\x_s)$, thus $x_u\in K := \pi(\K)$
and there are morphisms of racks
\begin{equation}\label{eq:reduction-unipotent}
\oc^{\K}_{\x_u}\simeq  \oc^{K}_{x_u} \hookrightarrow \Oc_x^{\Gb},
\end{equation}
the first by the isogeny argument and the second by
Remark \ref{rem:racks-typeD-flex} \ref{punto:flexD-b}.
Now the centralizer $C_{\Gsc^F}(\x_s)$ is a reductive subgroup of $\Gsc$ by \cite[Theorem 2.2]{Hu},
and $\K = C_{\Gsc}(\x_s)\cap \Gsc^F$.
Strictly, $C_{\Gsc^F}(\x_s)$ is not of Lie type in the sense above,
but close enough to allow some inductive
procedure. So, we are reduced to investigate the conjugacy classes
\begin{itemize}
  \item $x$ semisimple (the case $x = x_s$), or

  \item $x$ unipotent, and from this try to catch the general case.
\end{itemize}

\subsection{} In the first paper of the series, we deal with non-semisimple classes in $\Gb = \PSL_n(q)$, except $\PSL_2(q)$ with $q = 2,3,4,5,9$ which is either 
solvable or was treated in \cite{AFGV-ampa}; see \S \ref{subsec:intro-finite-Lie}.
To state our results, we start with some terminology. 
By the classical theory of the Jordan form, unipotent conjugacy classes in $\GL_n(q)$ are classified
by their type; $u \in \GL_n(q)$
is of type $\lambda = (\lambda_{1},\ldots, \lambda_{k})$
when the elementary
factors of its characteristic polynomial equal
$(X -1)^{\lambda_{1}}$, $(X -1)^{\lambda_{2}}$, \dots,
$(X -1)^{\lambda_{k}}$,
where $\lambda_{1}\geq \lambda_{2} \geq \cdots \geq \lambda_{k}$;
thus $u$ is  unipotent.

\newcounter{tabla}\stepcounter{tabla}
\renewcommand{\thetabla}{\Roman{tabla}}

\begin{theorem}\label{th:unipotent-slnq-collapse}

Let $x\in \Gb$ and pick $\x\in \SL_{n}(q)$ such that $\pi(\x) = x$, with Jordan
decomposition $\x = \x_s\x_u$.  Assume that $\x_u\neq e$. 
Then either  $\Oc = \oc^{\Gb}_{x}$ collapses or else $\x_s$ is central
and $\Oc$ is a unipotent class listed in  \ref{tab:uno}.
\end{theorem}
\begin{align}\label{tab:uno}\tag*{Table \thetabla}
\begin{tabular}{|c|c|c|c|}
\hline $n$  & type & $q$ & Remark  \\
\hline
$2$  & $(2)$ & even or not a square  & sober,  Lemma  \ref{lem:sl2} \\
\hline
$3$ &$(3)$  & 2  & sober,  Lemma  \ref{lem:sl2-ss-D} \eqref{step:sl2-tipo-tres}\\
\cline{2-4}
& $(2,1)$ & 2 &  cthulhu, Lemma  \ref{lem:sl2-ss-D} \eqref{step:sl2-tipo-dos-uno} \\
\cline{3-4}
&  & even $\ge 4$ & cthulhu, Prop. \ref{prop:2111-notD}, \ref{prop:21-notF}  \\
\hline
$4$ & $(2,1,1)$  & 2  & cthulhu, Lemma   \ref{lema:q=2-n=4}\\
\cline{3-4}
&  & even $\ge 4$ & not of type D, Prop. \ref{prop:2111-notD}, \\
& & & open for type F\\
\hline
\end{tabular}
\end{align}

\bigbreak  
Semisimple classes require different tools and are treated in work in progress.

\medbreak  
We deal with the Nichols algebras associated to the unipotent classes in \ref{tab:uno} in Lemma \ref{lem:sl2-YD-modules},
concluding the following result.

\begin{theorem}\label{theorem:YD-modules}
Let $\Oc$ be the conjugacy class of $x\in \Gb = \PSL_n(q)$ non-semisimple. Assume that either  $\Gb \neq \PSL_3(2)$, or else that $x$
is not of type $(3)$. Then $\dim\toba(\Oc,\rho)=\infty$, for every
$\rho\in \Irr C_{\Gb}(x)$.
\end{theorem}

The unipotent class $\Oc$ of type $(3)$ in  $\PSL_3(2)$ is sober and the centralizer of $x\in \Oc$  is cyclic of order 4.
Hence any abelian subrack has at most two elements. If $\rho\in \Irr C_{\Gb}(x)$ is given by $\rho(x) = -1$, then it is not possible
to decide whether the dimension of the Nichols algebra $\toba(\Oc, \rho)$ is finite or not by looking at subracks.

\medbreak  
Section \ref{sec:racks} is devoted to racks and Section \ref{sec:sln} to unipotent classes:
we prove Theorem \ref{th:unipotent-slnq-collapse} for them in \S \ref{subsec:summarize-collapse}.
In Section \ref{sec:non-ss} we prove the Theorem for non-semisimple classes, see Proposition \ref{prop:nonss-typeDPSLn}.

\subsection*{Notation} We denote the cardinal of a set $X$ by $\vert X\vert$. If $\ell$ is a positive integer, then we set 
$\I_{\ell} = \{i\in \N: 1\le i\leq \ell\}$.

Let $e_{i,j}\in \kk^{N\times P}$ be the  matrix with 1 in the position $(i,j)$ and 0 elsewhere.
We denote by $\id_N \in \kk^{N\times N}$ the identity matrix, and omit the subscript $N$ when clear from the context. 

Let $G$ be a group and $x_1, \dots, x_N \in G$. Then $\langle x_1, \dots, x_N\rangle$ denotes the subgroup generated by them.

\subsection*{Acknowledgements} We thank Istv\'an Heckenberger, who kindly communicated
us the proof of Theorem \ref{th:type4}, type F; Andrea Lucchini for information about conjugacy classes
and Leandro Vendramin, for many discussions on racks.
N. A. thanks Alberto de Sole for his hospitality during a visit to Universit\`a di Roma La Sapienza.

\section{Racks}\label{sec:racks}

\subsection{} A rack is a set $X \neq \emptyset$  with an
operation $\trid: X \times X \to X$ satisfying (a) $\varphi_x := x\trid \underline{\quad}$ is a
bijection for every $x \in X$,
and (b) the self-distributivity axiom $x\trid (y\trid z) =
(x\trid y) \trid (x\trid z)$ for all $x,y,z \in X$. Let $\Inn X$ be the subgroup of
$\Sim_X$ generated by $\varphi_x$, $x \in X$.
All racks in this paper are finite, unless explicitly stated.
The archetypical example of a rack is a conjugacy class $\Oc$
in a finite group $G$
with the operation $x\trid y = xyx^{-1}$, $x, y \in \Oc$.
We denote by  $\Oc_x^G$ (or $\Oc_x$ when no confusion arises)
the conjugacy class of $x$ in $G$.
Conjugacy classes are racks of a special sort, namely crossed sets,
as they satisfy
(c) $x\trid  x = x$ for all $x \in X$ and (d) $x\trid y = y$, iff $y\trid x = x$ for all $x,y \in X$,
see e.~g. \cite{AG-adv}. But this distinction
is not relevant for the purposes of this paper, so we assume that
\emph{all the racks appearing here are crossed sets}.
The following statement will be used along the paper.

\begin{obs}\label{rem:conjug-classes-normal-subgps}
Let $N$ be a normal subgroup of a finite group $G$, $x\in N$. Then there exists $x= x_1, \dots, x_s \in N$ such that
\begin{equation}\label{eq:conjug-classes-normal-subgps}
\Oc_x^G = \coprod_{1\le i \le s} \Oc_{x_i}^N,
\end{equation}
and $\Oc_{x}^N \simeq \Oc_{x_i}^N$ as racks for all $i \in \I_s$.
\end{obs}

Indeed, $\Oc_x^G \subset N$, since $N$ is normal, hence \eqref{eq:conjug-classes-normal-subgps} holds. Now, if $g_i\in G$ satisfies
$g_i\trid x = x_i$, then $g_i\trid \Oc_{x}^N = \Oc_{x_i}^N$ and the last claim follows.

\medbreak
A rack $X$ is \emph{abelian} when $x\trid y = y$, for all $x, y \in X$.
A rack is \emph{indecomposable} when it is not a disjoint union
of two proper subracks, or equivalently when it is a single $\Inn X$ orbit.
Any rack is the disjoint union of maximal indecomposable subracks (in a unique way),
called its indecomposable components \cite[1.17]{AG-adv}.

\medbreak
A rack $X$ is \emph{simple} when for any projection of racks
$\pi: X\to Y$, either $\pi$ is an isomorphism or
$Y$ has only one element. The classification of finite simple racks
is known \cite[3.9, 3.12]{AG-adv}, \cite{jo};
one of the main parts consists of conjugacy classes in a finite
simple non-abelian group.

\subsection{Racks of type D, F}\label{subsec:typeD}
We discuss criteria to decide that a rack  collapses, see Definition \ref{def:rack-collapses}.
We start by the relevant definitions. Let $G$ be a group and let $X$ be a finite rack.

\begin{definition}\label{def:rack-typeD}
 \cite[3.5]{AFGV-ampa} $X$ is \emph{of type D}
when it has  a decomposable subrack
$Y = R\coprod S$ with elements $r\in R$, $s\in S$ such that
\begin{equation}\label{eq:typeD-rack}
r\trid(s\trid(r\trid s)) \neq s.
\end{equation}
\end{definition}

\begin{obs}\label{rem:typeD-gps}
If $\Oc$ is a finite conjugacy class in $G$, then the
following are equivalent:

\begin{enumerate}
  \item\label{item:typeD-rack} The rack $\Oc$ is of type D.

  \item\label{item:typeD-group} There exist $r$, $s\in \Oc$
  such that $\Oc_r^{\langle r, s\rangle}\neq \Oc_s^{\langle r, s\rangle}$ and
\begin{equation}\label{eq:typeD-group}
(rs)^2\neq(sr)^2.
\end{equation}
\end{enumerate}
\end{obs}

\pf Notice that \eqref{eq:typeD-rack} and \eqref{eq:typeD-group}
are equivalent in this setting.
If \eqref{item:typeD-group} holds, then
$Y = \Oc_r^{\langle r, s\rangle}\coprod \Oc_s^{\langle r, s\rangle}$
is the desired decomposable subrack. Conversely if
\eqref{item:typeD-rack} holds with $Y = R\coprod S$ and $r\in R$, $s\in S$,
then $\Oc_r^{\langle r, s\rangle}\subset R$, $\Oc_s^{\langle r, s\rangle}\subset S$.
\epf

\begin{definition}\label{def:rack-typeF}
$X$ is \emph{of type F} if it has a family of subracks $(R_a)_{a \in \I_4}$
and a family $(r_a)_{a \in \I_4}$ with $r_a\in R_a$, and for $a\neq b \in \I_4$, $R_a \cap R_b = \emptyset$,
$R_a \triangleright R_b = R_b$,
\begin{align}
\label{eq:typeF-rack}
r_a\triangleright r_b &\neq r_b.
\end{align}
\end{definition}

Here F stands for a rack with a family of \emph{four} mutually
disjoint subracks.

\begin{obs}\label{rem:typeF-gps}
If $\Oc$ is a finite conjugacy class in $G$, then the
following are equivalent:

\begin{enumerate}
  \item\label{item:typeF-rack} The rack $\Oc$ is of type F.

  \item\label{item:typeF-group} There exist $r_a\in \Oc$, $a\in \I_4$,
  such that $\Oc_{r_a}^{\langle r_a: a\in \I_4\rangle}\neq \Oc_{r_b}^{\langle r_a: a\in \I_4\rangle}$, $a\neq b$ in $\I_4$, and
\begin{align}\label{eq:typeF-group}
r_a r_b &\neq r_br_a, & a & \neq b\in \I_4.
\end{align}
\end{enumerate}
\end{obs}

\pf Notice that \eqref{eq:typeF-rack} and \eqref{eq:typeF-group}
are equivalent in this setting.
If \eqref{item:typeF-group} holds, then
$R_a = \Oc_{r_a}^{\langle r_a: a\in \I_4\rangle}$, $a\in \I_4$
is the desired family of subracks. Conversely if
\eqref{item:typeF-rack} holds,
then $\Oc_{r_a}^{\langle r_a: a\in \I_4\rangle}\subset R_a$, for all $a\in \I_4$ and we have \eqref{item:typeF-group}.
\epf

\medbreak
The rack formulations \eqref{item:typeD-rack} in Remark \ref{rem:typeD-gps}, resp.  \ref{rem:typeF-gps}, are more
effective for applications to the classification of Hopf algebras, see Remark \ref{rem:racks-typeD-flex}; the 
equivalent formulations \eqref{item:typeF-group} are useful in proofs. 

\begin{obs}\label{rem:typeD-F-subgps}
Let $\Oc$ be a finite conjugacy class in $G$. If $\Oc$ is of type D, resp. F, then there is a maximal $K < G$
such that $\Oc\cap K$ is of type D, resp. F.
\end{obs}

Indeed, let  $r$, $s\in \Oc$  such that $\Oc_r^{\langle r, s\rangle}\neq \Oc_s^{\langle r, s\rangle}$; then
$\langle r, s\rangle \neq G$, so there is a maximal $K$ containing $\langle r, s\rangle$. Same for type F.

\medbreak
The following remark, a variation of \cite[Lemma 2.5]{FV}, is useful to check when 
the conditions in Remarks \ref{rem:typeD-gps} or  \ref{rem:typeF-gps} hold.

\begin{obs}\label{obs:rack-dihedral}
Let $G$ be a finite group and let $r, s\in G$ be involutions such that
$[s,\,r]\neq1$. Then $\Oc_r^{\langle r,\,s\rangle}\neq \Oc_s^{\langle r,\,s\rangle}$ if and only if $|rs|$ is even $>2$. 
\end{obs}

\begin{theorem}\label{th:type4}
A rack $X$ of type D (respectively, F) collapses.
\end{theorem}

\pf \textbf{Type D}: This is \cite[Theorem  3.6]{AFGV-ampa}, cf.
\cite[Theorem  8.6]{HS1}.

\medbreak
\textbf{Type F}: 
Let $q: X \times X \to \GL(n, \ku)$ be a finite faithful 2-cocycle on $X$.
We need to check that the Nichols algebra associated to the
braided vector space $(V, c) := (\ku X\otimes \ku^n, c^q)$ attached to $X$ and $q$
has infinite dimension. By hypothesis there is a subrack
$Y=\coprod_{a\in \I_4}R_a$ with $R_a \trid R_b = R_b$ as in Definition \ref{def:rack-typeF}.
Without loss of generality, we may assume that $Y = X$.
By \cite[6.14]{AG-adv}, cf. also \cite[Theorem 2.1]{AFGV-ampa},
$(V, c)$ can be realized as Yetter-Drinfeld module over a finite group $G$.
Actually we may choose the subgroup $G$ of $\GL(V)$ generated by
$g_x: V \to V$, $g_x(e_y\otimes v) = e_{x\trid y}\otimes q_{xy}(v)$, $x,y\in X$, $v\in \ku^n$.
Let $V_a := \ku R_a\otimes \ku^n$, a Yetter-Drinfeld submodule of $V$; clearly $V = \oplus_{a\in A} V_a$.
Now we may replace $V_a$ by a simple Yetter-Drinfeld submodule
$U_a$ with $r_a\in \supp U_a = \{g\in G: U_{a,g} \neq 0\}$, where $U_a = \oplus_{g\in G} U_{a,g}$ is the grading coming from the Yetter-Drinfeld
module structure.
Then $c^2 \neq \id$ on $U_a\otimes U_b$ for $a \neq b \in \I_4$ by \eqref{eq:typeF-rack}.
This means that the Weyl groupoid $\W$ of $U = \oplus_{a\in \I_4} U_a$, see \cite{AHS},
has rank at least 4 and the Dynkin diagram of one of his objects
would then have an edge between any two distinct vertices.
Now if $\dim \toba(X, q) < \infty$, then $\W$ is finite.
But  this contradicts the classification of finite Weyl groupoids
in \cite[Thm. 1.1]{CH}.
\epf

\medbreak

The proof for type F uses stronger facts than the proof for type D,
as it relies on the classification from \cite{CH}. 
By this reason, our order of preference for application of these criteria
is first type D, then F.

\begin{obs}\label{rem:racks-typeD-flex} Being open
conditions (i.~e., expressed by inequalities),
these notions enjoy some favorable properties.

\begin{punto}\label{punto:flexF-a} If a rack $X$ contains a
subrack of type D (respectively, F), then $X$ is
of type D (respectively,  F). If a rack $X$ projects onto a rack of type D (respectively,  F), then $X$ is
of type D (respectively, F).
\end{punto}

Let $K$ be a subgroup of $G$, $\tau\in K$, $C_G(K)$
the centralizer of $K$ in $G$.

\begin{punto}\label{punto:flexD-a} If $\Oc^K_\tau$ is of type D (respectively, F),
then so is $\Oc^G_\tau$. \end{punto}

\begin{punto}\label{punto:flexD-b} Let  $\kappa\in C_G(K)$.
The right multiplication by $\kappa$ identifies $\oc_\tau^K$
with a subrack of $\oc_{\tau\kappa}^G$;
if $\oc_\tau^K$ is of type D (respectively, F), then so is $\oc_{\tau\kappa}^G$.
 \end{punto}

\begin{punto}\label{punto:flexD-c} Assume that
$G = G_1 \times \dots \times G_r \ni x = (x_1, \dots, x_r)$. Then
$\Oc_x^G = \Oc_{x_1}^{G_1} \times \dots \times\Oc_{x_r}^{G_r}$;
hence, if $\Oc_{x_j}^{G_j}$ is of type D (respectively, F) for some $j$,
then so is $\Oc_{x}^{G}$.
 \end{punto}
 \end{obs}

Now an indecomposable rack $Z$ always admits a rack
epimorphism onto a simple rack $X$.
Therefore, any indecomposable rack having a quotient
simple rack of type D collapses.
Hence it is natural to ask for the classification of all \emph{simple} racks of type D or F.
See \cite{AFGV-simple} for the present status of this problem, in the case of type D.

\begin{lema}\label{lema:prod-racks-D}
Let $X$ and $Y$ be racks.

\medbreak
(i) Assume that there are $y_1\neq y_2\in Y$, $x_1\neq x_2\in X$
	such that $y_1\trid y_2 = y_2$, $x_1 \trid (x_2\trid (x_1\trid x_2)) \neq x_2$.
Then $X \times Y$ is of type D.

\medbreak
(ii)  Assume that there are $y_1,\dots, y_4\in Y$ all different, $x_1\dots, x_4 \in X$
  such that  $y_i\trid y_j = y_j$,  $x_i\trid x_j \neq x_j$ for $i\neq j\in \I_4$.
Then $X \times Y$ is of type F.

(iii) Let $X_i$ be disjoint sets provided with bijections $\varphi_i: X \to X_i$, $i\in \I_2$; 
$X^{(2)} := X_1 \coprod X_2 \simeq  X \times \I_2$ is a rack with  $\varphi_i(x) \trid \varphi_j(y) = \varphi_j(x \trid y)$, $i, j\in \I_2$.
If there are  $x_1\neq x_2\in X$ satisfying \eqref{eq:typeD-rack}, then $X^{(2)}$ is of type D.
\end{lema}

\pf Take $R = X\times \{y_1\}$, $S = X\times \{y_2\}$, $r = (x_1,y_1)$, $s = (x_2,y_2)$ in (i); 
$R_j = X\times \{y_j\}$, $r_j = (x_j, y_j)$, $j\in \I_4$, in (ii). Now (i) implies (iii).
\epf

\subsection{Cthulhu racks}\label{subsec:cthulhu}

Recall that a rack is \emph{cthulhu}
when it is neither  of type D nor of type F;
and that it is \emph{sober} if  every subrack is either abelian or indecomposable.
A sober rack is cthulhu.  More than this:

\begin{obs}\label{obs:rack-cthulhu-weak-sober} If all subracks \emph{generated by two elements} of a rack $X$
are either abelian or indecomposable, then $X$ is cthulhu.
\end{obs}

Here are some examples of these notions.

\begin{exa}\label{exa:rack-cthulhu-doble-tetraedro}
The  rack $\oc^{\sk}_3$ of 3-cycles in $\sk$, also known as
the cube rack, is the union of two tetrahedral racks (conjugacy
classes in $\ac$) not commuting with each other.
It is neither of type D  nor of type F.
\end{exa}

\begin{exa}\label{exa:rack-sober}
Every abelian rack  is sober. The tetrahedral rack is sober.
The conjugacy class of non-trivial unipotent elements in $\PSL_2(q)$,
 where either $q$ is even, or odd but not a square,
 is sober, cf. Lemma \ref{lem:sl2-cthulhu}.
\end{exa}

\begin{exa}\label{exa:rack-cthulhu}
The  rack of transpositions in $\sn$ is cthulhu for $n\geq2$
but not sober for $n\geq 4$; see \cite[Remark 4.2]{AFGV-ampa}
for other examples of conjugacy classes in
symmetric groups that are cthulhu.
\end{exa}

\begin{exa}\label{exa:rack-cthulhu-affine}
Let  $\D_n$  be the affine rack $(\Z_n, T)$ where $T$ is the
inversion;  when $n$ is odd, it is the class of involutions in the
dihedral group $\dn$ of order $2n$. If $n > 4$ is even, then
$\D_n$ is of type D \cite[Lemma 2.2]{AFGV-thr}. If $n$ is odd,
then $\D_n$ is sober. For, observe that every subgroup of  $\dn$
is either cyclic of order $d$ or isomorphic to a dihedral group
$\dl$, for some $d\vert n$. Let $X$ be a subrack of $\D_n$ and $H
=\langle X\rangle$. Since $X$ consists of involutions, $H \simeq
\dl$ for some $d\vert n$; hence $X$ is the class of involutions in
$H$, that is  indecomposable.

\end{exa}

\section{Unipotent classes in $\SL_n(q)$}\label{sec:sln}
Let $n\in \N$, $n\geq 2$. In this section, we consider $G = \SL_n(q)$
and investigate when a unipotent conjugacy class collapses.
By the isogeny argument,
the result carries over $\Gb = \PSL_{n}(q)$. We deal
with unipotent  classes of type D
in \S \ref{subsec:slnq-typeD},
with those of type F in \S  \ref{subsec:slnq-type4}.
We summarize in \S \ref{subsec:summarize-collapse}.

\medbreak
Before starting we state an observation useful not only in the unipotent context.
Let $u\in G$ with block decomposition
\begin{equation}\label{eq:u-unipot}
u =
\left( \begin{matrix}
         u_{1} & 0 &    \ldots  & 0 \\
         0 & u_{2} &  \ldots & 0 \\
       \vdots & & \ddots  & \vdots \\
                 0 &\ldots & \ldots & u_{k}
              \end{matrix}\right),
\end{equation}
where
$ u_{j} \in \SL_{\lambda_{j}}(q)$, $j \in \I_k$. By Remark \ref{rem:racks-typeD-flex}, we have:

\begin{lema}\label{lem:blockD}
If $ \Oc_{u_{i}}^{\SL_{\lambda_{i}}(q)} $
is of type D (respectively F) for some $i \in \I_k$, then so is $ \Oc_{u}^G$. \qed
\end{lema}

\subsection{Unipotent classes}\label{subsec:sln-unipotent}
Recall that a unipotent $u \in \GL_n(q)$
is of type $\lambda = (\lambda_{1},\ldots, \lambda_{k})$
when the elementary
factors of its characteristic polynomial equal
$(X -1)^{\lambda_{1}}$, $(X -1)^{\lambda_{2}}$, \dots,
$(X -1)^{\lambda_{k}}$,
where $\lambda_{1}\geq \lambda_{2} \geq \cdots \geq \lambda_{k}$.
A (unipotent)  $x\in \GL_n(q)$ (or its conjugacy class)
is \textit{regular} if it is of type
$(n)$, i.~e. if its characteristic and minimal polynomials coincide.
Every element of type
$\lambda = (\lambda_{1},\ldots, \lambda_{k})$ in $\GL_n(q)$
is conjugate to a $u$ with block decomposition as in \eqref{eq:u-unipot}
with $u_{i} = \left( \begin{matrix}
          1 & 1 &\ldots & 0 \\
       \vdots  & \ddots & \ddots & 0 \\
         0  & \ldots &  1 & 1 \\
         0  & \ldots & 0 & 1
              \end{matrix}\right)  \in \SL_{\lambda_{i}}(q)$. To describe unipotent conjugacy classes in
              $G = \SL_n(q)$ and other purposes we set some notation.
              For  $\ba = (a_{1},\ldots, a_{n-1}) \in\kc^{n-1}$,
define  $r_{\ba}$ and the set $R_{\ba} \subset G$ by:

\begin{align}\label{eq:def-Ra}
r_{\ba} = \left( \begin{matrix}
         1 & a_1 & 0 & \ldots  & 0 \\
         0 & 1 & a_2 &\ldots & 0 \\
       \vdots & & \ddots & \ddots & 0 \\
         0 &\ldots & \ldots &  1 & a_{n-1} \\
         0 &\ldots & \ldots & 0 & 1
              \end{matrix}\right) \in R_{\ba} = \left\{\left( \begin{matrix}
         1 & a_1& * & \ldots  & * \\
         0 & 1 & a_2 &\ldots & * \\
       \vdots & & \ddots & \ddots & * \\
         0 &\ldots & \ldots &  1 & a_{n-1} \\
         0 &\ldots & \ldots & 0 & 1
              \end{matrix}\right)  \right\}.
\end{align}
If $\ba = (a,1,\ldots, 1)$, $a\in \kc^\times$, then we simply write $r_{\ba} = r_{a}$.

\medbreak
The sets $R_{\ba}$ enjoy the following properties: $R_{\ba} R_{\bb}\subset R_{\ba + \bb}$, $R_{\ba}^{-1} = R_{-\ba}$, hence
$R_{\ba}\trid R_{\bb}  \subset R_{\bb}$. Thus, $\coprod_{\ba\in \Fc} R_{\ba}$ is a subrack of $G$ for every subset $\Fc$ of $\kc^{n-1}$.
We shall need more precise formulae. For $\ba = (a_{1},\ldots, a_{n-1})$, $\bb = (b_{1},\ldots, b_{n-1}) \in\kc^{n-1}$,
set
\begin{align}\label{eq:rack-ra-rb}
\theta^{k}_{\ba, \bb} &= a_{k}b_{k+1} - a_{k+1}b_{k},& &1\leq k \leq n-2,
\\
\gamma^{k}_{\ba, \bb} &= 2a_{k}b_{k+1} + (a_{k} + b_{k})(a_{k+1} + b_{k+1}), & &1\leq k \leq n-2,
\\ \label{eq:rack-ra-rb3}
\nu^{k}_{\ba, \bb} & = a_{k}b_{k+1}(a_{k+2}+ b_{k+2}) +  a_{k+1}b_{k+2}(a_{k}+ b_{k})& &1\leq k \leq n-3.
\end{align}

Then
\begin{align}\label{eq:rack-ra-trid}
r_{\ba}\trid r_{\bb} &=
\left( \begin{matrix}
         1 & b_1 & \theta^{1}_{\ba, \bb} &  -a_{3}\theta^{1}_{\ba, \bb}
& \ldots  & (-1)^{n-1}a_{3}\cdots a_{n-1}\theta^{1}_{\ba, \bb} \\
         0 & 1 & b_2 &\theta^{2}_{\ba, \bb} &
\ldots &  (-1)^{n-2}a_{4}\cdots a_{n-1}\theta^{2}_{\ba, \bb} \\
       \vdots & & \ddots & \ddots &  & \\
       0 &\ldots & \ldots &    & b_{n-2} &  \theta^{n-2}_{\ba, \bb}  \\
         0 &\ldots & \ldots &\ldots  & 1  & b_{n-1} \\
         0 &\ldots & \ldots & \ldots & 0 & 1
              \end{matrix}\right).
 \end{align}

Thus   $r_{\ba}\trid r_{\bb}\neq r_{\bb}$  if $\theta^{k}_{\ba,\bb}\neq 0$
for some $1\leq k \leq n-2$. Analogously,

\begin{align*}
(r_{\ba}r_{\bb})^2 &= \left( \begin{matrix}
         1 & 2(a_1+ b_1) & \gamma^{1}_{\ba, \bb} & \nu^{1}_{\ba, \bb} &  \ldots  & * \\
         0 & 1 & 2(a_2 + b_2) & \gamma^{2}_{\ba, \bb} &\nu^{2}_{\ba, \bb} & * \\
         0 & 0 & 1 & 2(a_3 + b_3) &\ddots & * \\
       \vdots & & & \ddots & \ddots & \gamma^{n-2}_{\ba, \bb} \\
         0 &\ldots & \ldots & 0 & 1 & 2(a_{n-1} + b_{n-1}) \\
         0 &\ldots & \ldots & 0  & 0 & 1
              \end{matrix}\right);
\end{align*}
hence $(r_{\ba}r_{\bb})^2 =(r_{\bb}r_{\ba})^2 $ implies that
\begin{equation}\label{eq:rarb^2}
\gamma^{k}_{\ba, \bb} = \gamma^{k}_{\bb, \ba} \iff 2\theta^{k}_{\ba,\bb}= 0,\
\forall\ 1\leq k \leq n-2 \text{ and }\newline
\nu^{k}_{\ba, \bb} = \nu^{k}_{\bb, \ba}\ \forall\ 1\leq k \leq n-3.
\end{equation}

\medbreak
Now every element in $G$ of type
$\lambda = (\lambda_{1},\ldots, \lambda_{k})$
is conjugate to one of the form \eqref{eq:u-unipot}
with $u_{i} = r_{{a}_i} \in
\SL_{\lambda_{i}}(q)$ for  some ${a}_i \in \F_q^\times$.

Indeed, assume that $V \in \SL_{n}(q)$ admits $C \in \GL_{n}(q)$ such that
$CVC^{-1}$ is of the form \eqref{eq:u-unipot} with regular unipotent blocks.
Consider  the diagonal matrix 
$D =  (\det C^{-1}, 1,\ldots, 1)\in (\F_q^\times)^{n}$.
Then $E = DC \in \SL_{n}(q)$ and $E V E^{-1}$ is of the form  \eqref{eq:u-unipot} with regular unipotent blocks.

\begin{obs}\label{obs:conj-class-unip-Fq} To study the unipotent conjugacy classes in $G$,
it suffices to consider  classes of elements of the form \eqref{eq:u-unipot}
with $u_{i} = r_{1}$, cf. Remark \ref{rem:conjug-classes-normal-subgps}.
\end{obs}

\medbreak
For further purposes, we shall need the following well-known description of the regular unipotent conjugacy classes in $G$.

\begin{lema}\label{lem:conj-class-unip-Fq}
Let $d := \GCF(q-1, n)$.
There are $d$ regular unipotent conjugacy classes in $G$, all isomorphic as racks.
Explicitly,
they are of the form $\oc_{r_{a}}$, for some $a\in \kc^\times$,
and $\oc_{r_{a}} = \oc_{r_{b}}$
if and only if $\theta^{n}a = b$ for some $\theta \in \F_{q}^{\times}$.
If $\ba = (a_{1},\ldots, a_{n-1})\in (\kc^{\times})^{n-1}$, then $R_{\ba} \subseteq \oc_{r_a}$ for $a =
a_{1}a_{2}^{2}\cdots a_{n-1}^{n-1}$.
\end{lema}

\pf
Let $x\in G$ be a regular unipotent element; we may assume that
\begin{equation}\label{eq:generic-unipotent-upptriang}
%\begin{align*}
x =
\left( \begin{matrix}
         1 & x_{12} & x_{13} &    \ldots  & x_{1n} \\
         0 & 1 & x_{23}  & \ldots & x_{2n} \\
       \vdots & & \ddots   & \ddots  & \vdots \\
         0 &\ldots & \ldots & 1 &x_{n,n-1} \\
  0 &\ldots & \ldots & 0 &1
             \end{matrix}\right),\   x_{i, i+1}\in \F_{q}^{\times}, \ 1\le i \le n-1.
%\end{align*}
\end{equation}
Let $a\in \F_{q}^{\times}$; we claim that $x \in \oc_{r_a}$ if and only if
 \begin{align}
 \label{eq:conj-class-lambda}
\theta^{n} a &= x_{12}x_{23}^{2}x_{34}^{3}\cdots x_{n-1,n}^{n-1}&
\text{ for some }&\theta \in \F_{q}^{\times}.
 \end{align}
Indeed, $x \in \oc_{r_{a}}$ if and only if
there exists $C = (c_{ij}) \in \SL_n(q)$ such that
$C r_a = x C$ which holds if and only if the following
linear equations hold
\begin{eqnarray}
&c_{n,j}  = 0 &\text{for all }1\leq j<n,\label{eq:conj1}\\
&c_{i,j}  = \sum_{k=i+1}^{n} x_{i,k}c_{k,j+1} &\text{for all } 1\leq i<n, 2\leq j <n,
\label{eq:conj2}\\
&a c_{i,1}  = \sum_{k=i+1}^{n} x_{i,k}c_{k,2} &\text{for all } 1\leq i<n,
\label{eq:conj3}\\
&0  = \sum_{k=i+1}^{n} x_{i,k}c_{k,1} &\text{for all } 1\leq i<n. \label{eq:conj4}
\end{eqnarray}
By a direct computation using \eqref{eq:conj1}, \eqref{eq:conj2}
and \eqref{eq:conj4}, $c_{ij} = 0$ for all $1\leq j < i \leq n$,
i.~e. $C$ is upper triangular. Thus, $a c_{11} = x_{12}c_{22}$
from \eqref{eq:conj3}, and $c_{ii} = x_{i,i+1}c_{i+1,i+1}$
for all $1< i <n$, from \eqref{eq:conj2}. Since $\det C = 1$, 
\begin{equation}\label{eq:conj-lambda-x-a}
a = a c_{11}\cdots c_{nn} =
x_{12}x_{23}^{2}\cdots x_{n-1,n}^{n-1}c_{n,n}^{n}.
\end{equation}
Thus, if $x\in \oc_{r_a}$, then it is conjugated to $r_a$ by an upper triangular matrix
$C$ and  \eqref{eq:conj-class-lambda} holds with $\theta = c_{n,n}^{-1}$.
Conversely, if \eqref{eq:conj-class-lambda} is satisfied, then
define an upper triangular matrix $C$ by
$c_{n,n} = \theta^{-1}$, $c_{ii} = x_{i,i+1}c_{i+1,i+1}$
for $1< i <n$, $c_{11}= a^{-1}x_{12}c_{2,2}$ and use
equations \eqref{eq:conj2} to find the remaining elements.
Consequently,
$\oc_{r_{a}} = \oc_{r_{b}}$
if and only if $\theta^{n}a = b$
for some $\theta \in \F_{q}^{\times}$; i.~e. the set of regular unipotent  classes in $G$
is parameterized by the quotient of the cyclic
group $\F_{q}^{\times}$ by the image of the map
by $x \mapsto x^n$. Since the kernel of this map has order $d=\GCF(n,q-1)$,
we get $d$ different classes.
\epf

\subsection{Unipotent conjugacy classes in $\PSL_2(q)$}\label{subsec:sln2-aparte}
We start with unipotent classes in $\PSL_2(q)$; here $q\neq 2$, $3$, $4$, $5$, $9$, see \S  \ref{subsec:intro-finite-Lie}.
First we recall Dickson's classification of all subgroups of $\PSL_2(q)$.
Let $d = (2, q-1)$.

\begin{theorem}\label{th:subgps-psl2}\cite[Theorems 6.25, p. 412; 6.26, p. 414]{suzuki}
A  subgroup of $\PSL_2(q)$ is isomorphic to one of the following groups.

\begin{enumerate}\renewcommand{\theenumi}{\alph{enumi}}\renewcommand{\labelenumi}{(\theenumi)}
  \item\label{dickson:point-dihedral} The dihedral groups of order $2(q\pm 1)/d$ and their subgroups. There are always such subgroups.
  \item\label{dickson:point-elementary-abelian}
  A group $H$ of order $q(q-1)/d$ and its subgroups. It has a normal $p$-Sylow subgroup $Q$ that is elementary
  abelian and the quotient $H/Q$ is cyclic of order $(q-1)/d$. There are always such subgroups.
  \item\label{dickson:point-a4} $\mathbb A_4$, and there are such subgroups except when $p=2$ and $m$ is odd.
    \item\label{dickson:point-s4} $\mathbb S_4$, and there are such subgroups if and only if $q^2 \equiv 1 \mod 16$.
      \item\label{dickson:point-a5} $\mathbb A_5$, and
      there are such subgroups if and only if $q(q^2 - 1) \equiv 0 \mod 5$.
  \item\label{dickson:point-psl2} $\PSL_{2}(t)$ for some
  $t$ such that $q = t^h$, $h\in \N$.
  There are always such subgroups.
  \item\label{dickson:point-pgl2} $\PGL_{2}(t)$ for some
  $t$ such that $q = t^h$, $h\in \N$. If $q$ is odd, then
 there are such subgroups  if and only if $h$ is even and $q = t^h$. \qed
\end{enumerate}
\end{theorem}

\begin{lema}\label{lem:sl2-cthulhu}
Assume that $q$ is either even, or else odd but not a
square. Then a unipotent conjugacy class $\oc$ of $\PSL_2(q)$
is sober, hence cthulhu.
\end{lema}

\pf Let $X$ be a subrack of $\oc$; we show that $X$ is
either abelian or indecomposable.
Let $K$ be the subgroup of $\PSL_2(q)$ generated by $X$.
Since $X$ generates $K$, it is a union of (unipotent)
$K$-conjugacy classes \cite[1.9]{AG-adv}. We may
assume that $r_1\in X$.
The order of any element in $X$ is $p$,
so $p$ divides $\vert K\vert$; this excludes case
\eqref{dickson:point-dihedral} in Theorem
\ref{th:subgps-psl2} for $p$ odd.
Assume that $q$ is even, so that $d=1$, and $K$
is a dihedral group
of order $2(q\pm 1)$. Then $X$ is the rack of
involutions of $K$, which is indecomposable, see
Example \ref{exa:rack-cthulhu-affine}.

If $K$ is as in case \eqref{dickson:point-elementary-abelian},
then $X \subset Q$, hence it is an abelian rack.
If $K\simeq \mathbb A_4 \simeq \PSL_2(3)$,
case \eqref{dickson:point-a4}, then $p=2$ or $3$.
If $p=2$, then $X$ could not generate $K$,
being contained in the normal 2-Sylow subgroup of $K$.
If $p = 3$, then we are reduced to case
\eqref{dickson:point-psl2}.
If $K\simeq \mathbb S_4$, case \eqref{dickson:point-s4},
then $p=2$ or $3$; but the $3$-cycles in $\sk$
generate $\ac$ so $3$ is not possible, whereas $p=2$
is excluded by Theorem \ref{th:subgps-psl2}.
If $K\simeq \mathbb A_5 \simeq \PSL_2(5)$, case
\eqref{dickson:point-a5}, then $p=2$, $3$ or $5$.
If $p=2$, then $X$ is indecomposable, being the unique
class of involutions in $\aco$.
If $ p = 3$, then $p^{2m} - 1 \equiv 0 \mod 5 \iff m$ is
even, excluded by hypothesis.
If $p = 5$, then we are reduced to case \eqref{dickson:point-psl2}.

Assume then that  $K \simeq \PSL_{2}(t)$, $q =t^h$,
case \eqref{dickson:point-psl2}.
For $q$ even,
$\PGL_{2}(t) \simeq \PSL_{2}(t)$ has just one regular
unipotent conjugacy class, so $X$ is indecomposable
by \cite[1.9, 1.15]{AG-adv}; for $\PSL_{2}(2)\simeq  \st$ this is clear.
Assume that $q$ is odd.
Let $s \in X$; is $K$-conjugate to  $r_x$ for some $x\in\F_{t}^\times$.
But $\Oc_s^{G} = \Oc = \Oc_{r_1}^{G}$, hence $x \in (\F_q^\times)^2$.
Since $m$ is odd,
this only happens when $x$ is a  square in $\F_{t}^\times$,
i. e. when $\Oc_s^{K}=\Oc_{r_1}^{K}$.
Hence $X = \Oc^{K}_{r_1}$ is indecomposable.
Here we have to argue separately for
 $\PSL_{2}(3)\simeq \ac$, but in this case the claim is clear.

Finally, case \eqref{dickson:point-pgl2} is excluded
when $q$ is odd because $q$ is not a square.  \epf

\begin{lema}\label{lem:sl2} A non-trivial unipotent conjugacy class
$\Oc$ in $G = \SL_{2}(q)$, respectively in  $\PSL_{2}(q)$, is of type D if and only if $q > 9$ is an odd square.
\end{lema}

We excluded $\PSL_{2}(9) \simeq \mathbb{A}_6$ but in this case $\oc_{r}$ is
not of type D by \cite[Remark 4.2 (b)]{AFGV-ampa}.

\pf By Remark \ref{obs:conj-class-unip-Fq},
we may assume that
$\Oc = \Oc_r$ with $r= r_{1}$.
Suppose  $q \neq 9$ is an odd square.
Let $x\in \F_p^\times - (\F_p^\times)^{2}$; since $q\neq9$
we may assume that $x\neq2$. Let $v= r_x$.
         Since $x$ is not a square in $\F_p^\times$,  $\Oc_{r}^{G} =
\Oc_{v}^{G}$ but $R:=\Oc_{r}^{\SL_2(p)} \neq
\Oc_{v}^{\SL_2(p)}=:S$. Let
$s= \left(\begin{smallmatrix}
          0& 1\\
	  -1&0
         \end{smallmatrix}\right)
         \left(\begin{smallmatrix}
          1& x\\
	  0&1
         \end{smallmatrix}\right)
         \left(\begin{smallmatrix}
          0& -1\\
	  1&0
         \end{smallmatrix}\right)=\left(\begin{smallmatrix}
          1& 0\\
	  -x&1
         \end{smallmatrix}\right)\in S.$
         Then $(rs)^{2}\neq (sr)^{2}$ showing that $\Oc^{G}_{r}$ is of type D.
Conversely, if $q$ is not an odd square, then $\Oc^{G}_r$ is cthulhu
by Lemma \ref{lem:sl2-cthulhu}, hence not of type D.
\epf

The exceptional isomorphism $\PSL_3(2) \simeq \PSL_2(7)$ motivates the analysis of some semisimple classes in this last group.

\begin{lema}\label{lem:sl2-ss-D}
Let $\Oc$ be the conjugacy class of $x\in \PSL_2(7)$.

\begin{enumerate}\renewcommand{\theenumi}{\alph{enumi}}\renewcommand{\labelenumi}{(\theenumi)}
	
	\item \label{step:sl2-tipo-dos-uno}
	If $\ord x = 2$, then $\Oc$ is  cthulhu.
	
	\item\label{step:sl2-tipo-tres}
If $\ord x = 4$, then $\Oc$  is sober.
	
\end{enumerate}

Hence the class of  
type $(2,1)$, respectively $(3)$, in $\PSL_3(2)$ is cthulhu, respectively sober.
\end{lema}

\pf The proper subgroups of $\PSL_2(7)$  are isomorphic either to  $\mathbb D_3\simeq \st$,  $\mathbb D_4$, 
the non-abelian group of  order 21, $\ac$, or $\sk$, or their subgroups.
Let $X$ be a subrack of $\oc$ and $K  = \langle X\rangle$; we show by inspection that $X$ is either abelian, or indecomposable,
or the union of at most 3 subracks that do not fulfill \eqref{eq:typeD-rack}.
Suppose that $\ord x = 2$. First, $\oc_{(12)}^{\st}$ is indecomposable. Second, $\mathbb D_4 = \langle r,s\vert r^4 =s ^2 = \id, srs =r^3\rangle$
has 3 classes of involutions: $\{r^2\}$ which is central, and the abelian racks $\{s, sr^2\}$, $\{sr, sr^3\}$ not commuting with each other;
\eqref{eq:typeD-rack} does not hold here.
The involutions of $\ac$ generate the 2-Sylow subgroup, so $K = \ac$ is not possible.
If $K = \sk$, then
$\oc_{(12)(34)}^{\ac} = \oc_{(12)(34)}^{\sk}$ does not generate $K$,  $\oc_{(12)}^{\sk}$ is indecomposable by \cite[3.2 (2)]{AG-adv},
and \eqref{eq:typeD-rack} does not hold in $\oc_{(12)(34)}^{\sk}\coprod \oc_{(12)}^{\sk}$. 
So, $\oc$ is cthulhu but not sober.
If $\ord x = 4$, then $K$ could be either abelian, $\mathbb D_4$ or $\sk$.
The elements of order 4 in $\mathbb D_4$ generate a cyclic subgroup.
If $K \simeq \sk$,  then $X = \oc^{\sk}_4$ is indecomposable.
\epf

\subsection{Unipotent classes of type D} \label{subsec:slnq-typeD}

\subsubsection{Odd characteristic} \label{subsubsec:slnq-typeD-odd}
Here we assume that  $q$ is odd.

\begin{lema}\label{cor:rack-typeD-SLn}
Let $u \in G$ be a unipotent element of
type $(\lambda_{1},\ldots ,\lambda_{k})$.
If $\lambda_{1} >2$, then
the conjugacy class $\oc_{u}$ is of type D.
\end{lema}

\pf Assume that $u$ is regular, i.~e. $\lambda_{1} = n > 2$.
By Remark \ref{obs:conj-class-unip-Fq}
we may suppose that  $\Oc = \oc_{r_{\uno}}$.
Let $\zeta\in\kc^\times$ with $\zeta^3\neq1$,
$t$ the diagonal matrix $(1,\zeta,\zeta^{-1},1,\ldots,1)$ and
$\bb=(\zeta^{-1},\zeta^2,1,\ldots,1)$.
Then $tr_{\uno}t^{-1}\in R_{\bb}$ and $R_{\uno}\coprod R_{\bb}$
is a decomposable subrack of $\oc_{r_{\uno}}$.
Besides,  $r_{\uno}\trid (r_{\bb}\trid (r_{\uno}\trid r_{\bb})) \neq r_{\bb}$
by \eqref{eq:rarb^2}, so that $\oc_{r_{\uno}}$ is of type D.
In the general case, we may assume that $u$ is as in \eqref{eq:u-unipot} with
$ u_{j} \in \SL_{\lambda_{j}}(q)$, $j  \in \I_k$.
Then  Lemma \ref{lem:blockD} applies.
\epf

We next deal with non-trivial unipotent conjugacy classes not covered by
the previous lemma, that is those of
type $(2, 2,\ldots, 1, \ldots , 1)$.

\begin{lema}\label{lem:2-1,2-2}
Let $u \in G$ be a unipotent element. Assume that either
\renewcommand{\theenumi}{\alph{enumi}}\renewcommand{\labelenumi}{(\theenumi)}
\begin{enumerate}
\item\label{item:2-2}  $n =4$ and $u$ has type $(2,2)$ or
\item\label{item:2-1} $n =3$ and $u$ has type $(2,1)$.
\end{enumerate}
Then the conjugacy class $\oc_{u}$ is of type D.
\end{lema}

\pf
\eqref{item:2-2}: We may assume
$u= \left(\begin{smallmatrix}
     r_{1} & 0\\
     0 & r_{1}
    \end{smallmatrix}  \right)$, Remark \ref{obs:conj-class-unip-Fq}. Let $\zeta\in\F_q^\times - (\F_q^\times)^2$,
$r=u$, $t= \left(\begin{smallmatrix}
     0 & -\zeta & 1 &1\\
     1 & 0 & 0 & 0\\
     0  & 0 & 0 & -\zeta^{-1}\\
     0 & 0 & 1 & 0\\
    \end{smallmatrix}  \right)$ and
$s= t\trid u = \left(\begin{smallmatrix}
     1 & 0 & -\zeta  & 0 \\
     -\zeta^{-1} & 1 & -1 & \zeta^{-1} \\
     0  & 0 & 1 & 0\\
     0 & 0 & -\zeta & 1\\
    \end{smallmatrix}  \right)$.
A direct computation shows that $(rs)^2\neq(sr)^2$. Moreover,
$\langle r,\,s\rangle$ is strictly contained in the subgroup  $H$ of $\SL_4 (q)$
of block upper triangular matrices with diagonal blocks in $\SL_2 (q)$.
Since
$\left(\begin{smallmatrix}
     1 & 0 \\
     -\zeta^{-1} & 1
    \end{smallmatrix}  \right)$ is conjugated
to
$\left(\begin{smallmatrix}
     1 & \zeta^{-1} \\
     0 & 1
    \end{smallmatrix}  \right)$
in $\SL_{2}(q)$ and this is not conjugated to $r_{1}$
 by Lemma \ref{lem:conj-class-unip-Fq}, it turns out that
$\Oc_r^H\neq\Oc_s^H$, so $\Oc_r^{\SL_4 (q)}$ is of type D.

\eqref{item:2-1}: We may assume that
$u= \left(\begin{smallmatrix}
     1 & 1 & 0\\
     0 & 1 & 0\\
     0 & 0 & 1
    \end{smallmatrix}  \right)$.
Take $r=u$ and $s=\left(\begin{smallmatrix}
     0 & 1 & 0\\
     0 & 0 & 1\\
     1 & 0 & 0
    \end{smallmatrix}  \right) \trid r = \left(\begin{smallmatrix}
     1 & 0 & 0\\
     0 & 1 & 0\\
     1 & 0 & 1
    \end{smallmatrix}  \right) \in \Oc_r^{\SL_3 (q)}$.
Then $(rs)^2\neq(sr)^2$, since $q$ is odd.
Also, $H=\langle r,\, s\rangle=\left\{   \left(\begin{smallmatrix}
     1 & a & 0\\
     0 & 1 & 0\\
     c & b+ac & 1
    \end{smallmatrix}  \right)~|~a,b,c\in\F_p\right\}$
because $[r^k,s^m] =  \left(\begin{smallmatrix}
     1 & 0 & 0\\
     0 & 1 & 0\\
     0 & km & 1
    \end{smallmatrix}  \right)$
and $[r,[r^k,s^m]]=[s,[r^k,s^m]]=1$. It is not hard to verify that $\Oc_r^H\neq\Oc_s^H$.
\epf

If $\oc_{u}^{\SL_n(q)}$ is of type D, then
$\oc_{u}^{\GL_n(q)}$ is of type D too. Thus,
the previous results
apply to nontrivial unipotent conjugacy classes in $\GL_n(q)$
with the prescribed hypothesis.  We deal with the remaining cases.

\begin{lema}\label{lem:GL(2,q)}
If $u= \left(\begin{smallmatrix} 1 & 1\\ 0 & 1 \end{smallmatrix}\right)$
and $q> {3}$,
then $\oc_{u}^{\GL_2(q)}$ is of type D.
\end{lema}
\pf
Consider the subsets of
$\oc_{u}$ given by
\begin{align*}
R & = \{AuA^{-1}:\ A\in \GL_2(q), \det A \not\in(\F_{q}^{\times})^{2}\},\\
S & = \{AuA^{-1}:\ A\in \GL_2(q), \det A \in(\F_{q}^{\times})^{2}\}.
\end{align*}
Let $A, B\in \GL_2(q)$ and set $r = AuA^{-1}$ and $s=BuB^{-1} \in \oc_{u}$. Then
$$r\trid s= (AuA^{-1})(BuB^{-1})(AuA^{-1})^{-1} =
(AuA^{-1}B)u(B^{-1}Au^{-1}A^{-1})$$
and $\det (AuA^{-1}B) =
\det B$. Hence $R, S$ are subracks of $\oc_{u}$,
$R\trid S \subseteq S$ and $S\trid R\subseteq S$. Moreover,
$R\cap S \neq \emptyset$ if and only if
there exists $B\in C_{\GL_2(q)}(u)$ with $\det B $ not a square.
But $C_{\GL_2(q)}(u)=\{\left(\begin{smallmatrix}
                         a & b\\ 0 & a
                        \end{smallmatrix}\right):\ a,b\in \F_{q} \}$, so that
both racks are disjoint and $Y = R\coprod S$ is a decomposable
subrack of $\oc_{u}$.
Now let $d\in \F_{q}$ be not a square and $r_{d} = d^{-1}\left(\begin{smallmatrix}
                         d-1 & 1\\ -1 & d+1
                        \end{smallmatrix}\right) = A_{d}u A_{d}^{-1} \in R$ with
$A_{d} = \left(\begin{smallmatrix}
                         1 & 0\\ 1 & d
                        \end{smallmatrix}\right)$.
If we set $ s=u $, a direct computation in $ \GL_2(q) $ shows
that
$$ s\trid (r_{d}\trid s) =  d^{-2}\left(\begin{smallmatrix}
                         d^{2} + d-2 & d^{2}-4d + 4\\ -1 & d^{2}-d+2
                        \end{smallmatrix}\right)=t.$$
Hence, $r_{d}\trid(s\trid (r_{d}\trid s)) = s$ if and only if
$ t = r_{d}^{-1}\trid s$ if and only if
$$ \left(\begin{smallmatrix}
                         d^{2} - d-1 & (d+ 1)^{2}\\ -1 & d^{2}+ d +1
                        \end{smallmatrix}\right) =  \left(\begin{smallmatrix}
                         d^{2} + d-2 & (d-2)^{2}\\ -1 & d^{2}-d+2
                        \end{smallmatrix}\right),$$
if and only if $ (d+1)^{2} = (d-2)^{2} $ and
$ 2d = 1 $ in $ \F_{q} $. Thus, $r_{d}\trid(s\trid (r_{d}\trid s)) \neq s$ if
$ 2d\neq 1 $. If $q\neq 3$, such a $ d $ always exists,
showing that $ \oc_{u} $ is of type $ D $.
\epf

\subsubsection{Even characteristic} \label{subsubsec:slnq-typeD-even}
Here we assume that  $q$ is even.

\begin{lema}\label{lem:difblocksallq} 
Let $u \in G$ be a unipotent element of type $\lambda=(\lambda_{1}, \ldots, \lambda_{k})$; 
assume that  $\lambda_{i} \geq \lambda_{i+1} \geq 3$ for some $1\leq i\leq k-1$.
Then the conjugacy class $\oc := \oc_{u}$ is of type D.
\end{lema}

\pf By Lemma \ref{lem:blockD}, it is enough to look at the following specific unipotent class:
If $\lambda = (\lambda_{1}, \lambda_{2})$ with $\lambda_{1}\geq \lambda_{2}\geq 3$, then $\Oc$ is of type D.

Let $x_{i} = r_{(1,\ldots, 1)}\in \F^{\lambda_{i}\times\lambda_{i}}$, $i= 1,2$. 
By Remark \ref{obs:conj-class-unip-Fq} we may assume
that  $u=\left(\begin{smallmatrix}x_{1}&0\\
0&x_{2}\end{smallmatrix}\right)$.
Notice that $x_1^{-1} =\left(\begin{smallmatrix} 1 & 1 &    \ldots  & 1 \\
         0 & 1 &  \ldots & 1 \\
       \vdots & & \ddots  & \vdots \\
                 0 &\ldots & \ldots & 1\\
\end{smallmatrix}\right)\in \Oc_{x_1}$ by Lemma \ref{lem:conj-class-unip-Fq}.
Let
\begin{align*}
R_1 &= \left\{\left(\begin{smallmatrix} x_1& Z\\
0& x_2\\\end{smallmatrix}\right) ~|~ Z = (z_{ij})\in \kdos^{\lambda_1\times \lambda_2} \right\},& R &= R_1 \cap \oc; \\
S_1 &= \left\{\left(\begin{smallmatrix} x_1^{-1} & Z\\
0 & x_2\\\end{smallmatrix}\right) ~|~ Z = (z_{ij})\in \kdos^{\lambda_1\times \lambda_2} \right\},& S &= S_1 \cap \oc.
\end{align*}
Since $Y_1 = R_1 \coprod S_1$ is a decomposable subrack, so is $Y = R\coprod S \subset\oc$. Let
\begin{align*}
z_1 &= \left(\begin{smallmatrix} 1& 0 &0 &\dots & 0\\ 1& 0 &0 &\dots & 0\\ \vdots& \vdots &\vdots & & \vdots\\
              1& 0 &0 &\dots & 0 \\0 & 1& 0  &\dots & 0
             \end{smallmatrix}\right),& 
             z_2 &= \left(\begin{smallmatrix} 0& 0 &0 &0 &\dots & 0\\ 0& 0 &0 & 0&\dots & 0\\ \vdots& \vdots &\vdots &\vdots & & \vdots\\
             0& 0 &1 &1 &\dots & 0 \\0& 0& 0 & 1  &\dots & 0
             \end{smallmatrix}\right),&
            z_3 &=  \left(\begin{smallmatrix} 0& 0 &1 &1 &0 &\dots & 0\\ 0& 0 &1 & 1& 0& \dots & 0\\ \vdots& \vdots &\vdots &\vdots &\vdots & & \vdots\\
             0& 0 &1 &1 &0 &\dots & 0 \\0& 0& 0 & 1  & 1& \dots & 0
             \end{smallmatrix}\right).
\end{align*}
If $P = \left(\begin{smallmatrix} \id_{\lambda_1} & e_{\lambda_1,1}\\
0&\id_{\lambda_2} \end{smallmatrix}\right) = P^{-1} \in G$, then $v = P\left(\begin{smallmatrix} x_1^{-1}&0\\
0&x_2 \end{smallmatrix}\right)P = \left(\begin{smallmatrix} x_1^{-1}& z_{1}\\
0&x_2 \end{smallmatrix}\right)\in S$,
\begin{align*}
(uv)^2 &= \left(\begin{smallmatrix} \id & x_1 z_{1}\\ 0&x_2^2\end{smallmatrix}\right)^2
= \left(\begin{smallmatrix} \id & x_1 z_{1} (\id +  x_2^2)\\ 0&x_2^4\end{smallmatrix}\right)=
\left(\begin{smallmatrix} \id & z_{2} \\ 0&x_2^4\end{smallmatrix}\right), \text{ and }\\
(vu)^2 &= \left(\begin{smallmatrix} \id & z_{1}x_2 \\ 0&x_2^2\end{smallmatrix}\right)^2 =
\left(\begin{smallmatrix} \id & z_{1}x_2 (\id +  x_2^2)\\ 0&x_2^4\end{smallmatrix}\right) =
\left(\begin{smallmatrix} \id & z_{3}\\ 0&x_2^4\end{smallmatrix}\right)\text{ when } \lambda_2 > 3; \text{ but  }
\end{align*}
when $\lambda_2 = 3$, $(uv)^2 = \left(\begin{smallmatrix} \id & e_{2,3} \\ 0&x_2^4\end{smallmatrix}\right)$, 
$(vu)^2 = \left(\begin{smallmatrix} \id & e_{1,3} + e_{2,3} \\ 0&x_2^4\end{smallmatrix}\right)$. Thus $\Oc$ is of type D.

\epf

\begin{lema}\label{lema:q=2-n=4}
Let $u \in G$  of
type $(\lambda_{1},\ldots ,\lambda_{k})$ and assume that either
\begin{enumerate}\renewcommand{\theenumi}{\alph{enumi}}\renewcommand{\labelenumi}{(\theenumi)}
	
	\item\label{item:typeD-even-cuatro}  $\lambda_{1} = 4$.

	\item\label{item:typeD-even-tres-uno}  $\lambda_{1} = 3$, $\lambda_{2} = 1$.
	
	\item\label{item:typeD-dos-dos} $\lambda_{1} =\lambda_{2} = 2$.
\end{enumerate}

Then the conjugacy class $\oc_{u}$ is of type D. Furthermore, the class in $\PSL_4(2)$ of type $(2,1,1)$ is cthulhu.
\end{lema}

\pf Assume that $n = 4$, $q = 2$: Since $\PSL_4(2) \simeq \ao$,
we apply \cite{AFGV-ampa}.
There are two classes of involutions in $\ao$, of
types $(1^4, 2^2)$ or $(2^4)$; with centralizers of orders 96 and 192, respectively.
The former is of type D \cite[Table 2]{AFGV-ampa}, and
the latter is cthulhu because its proper subracks generated
by two elements are abelian racks and dihedral racks with 3
and 4 elements
\cite[4.2 (f)]{AFGV-ampa}. Now the class in $\PSL_4(2)$
of type $(2,1,1)$, respectively type $(2,2)$, has centralizer of order 192 so it is cthulhu,
respectively of order 96 and so is of type D.
Also, there are two classes of elements of order
4 in $\ao$, of types $(1^2, 2, 4)$ or $(4^2)$, 
both of type D \cite[Table 1 and Step 9]{AFGV-ampa}.
Hence the classes in $\PSL_4(2)$ of types $(3,1)$ and $(4)$ are of type D.
Now the claim for the classes $(2,2)$, $(3,1)$ and $(4)$, for $q$ even, follows
as $\SL_{n}(q)<\SL_{n}(q^{j})$ for any $j\in \N$; here Remark \ref{obs:conj-class-unip-Fq} is needed.
Finally Lemma \ref{lem:blockD} applies.
\epf

We now present a negative result.

\begin{prop}\label{prop:2111-notD} The unipotent classes of type $(2,1,1,1,1,\dots)$ in $\SL_n(q)$  for $q$ even and $n\geq 2$ 
are not of type D.\end{prop}
\pf Let $\Oc$ be a  class of type $(2,1,1,1,1,\dots)$ in $G = \SL_n(q)$. Let  $\U^F$,
respectively $\T^F$, the subgroup of unipotent upper-triangular, respectively diagonal matrices, in $G$.
Without loss of generality we may assume that it is represented by 
$r=\Id_n+e_{1,n}$, which lies in $Z(\U^F)$. We will show that if  $s\in \Oc$ satisfies 
$[s,\,r]\neq1$ and $\Oc_r^{\langle r,\,s\rangle}\neq \Oc_s^{\langle r,\,s\rangle}$, then $(rs)^2=(sr)^2$. 
Let $s\in\Oc$ satisfy $[r,\,s]  \neq 1$, and let $g\in G$ such that
$s=grg^{-1}$. By  \cite[24.1]{MT} $g$ can be decomposed as $g=u n_w tv$ where $n_w$ is a monomial matrix with 
coefficients in $\F_2$; $u,\,v\in \U^F$ and $t\in\T^F$. Then $s=u n_w t r t^{-1}n_w^{-1}u^{-1}$. We have:
\begin{align*}
t r t^{-1} &= \Id_n+\xi e_{1,n},& &\mbox{ for some } \xi\in\F_q^{\times},\\
\sigma :&=n_w t r t^{-1}n_w^{-1}=\Id_n+\xi e_{i,j},& &\mbox{ for some } i\neq j.
\end{align*}
Further, $[r,\,s]\neq1$ iff $u^{-1}[r,\,s]u=[r,\, \sigma]\neq1$ and this happens only if  either $i=n$ or $j=1$, or both. 
Assume first $(i,\,j)=(n,\,1)$. Since $r\in Z(\U^F)$ we have
$K:=\langle r,\,s\rangle\simeq u^{-1} \langle r,\, s\rangle u=\langle r,\,\sigma\rangle\simeq \langle \left(\begin{smallmatrix}
1&1\\
0&1
\end{smallmatrix}\right),\,\left(\begin{smallmatrix}
1&0\\
\xi&1
\end{smallmatrix}\right)\rangle\leq \SL_2(q)$. By Lemma \ref{lem:sl2-cthulhu},
$\Oc_r^K=\Oc_s^K$. 
Assume now $i=n$ and  $j\neq 1,\,n$.  
Then 
$$\begin{array}{ll}
u^{-1}(rs)^2u&=(r\sigma)^2=((\Id_n+e_{1,n})(\Id_n+\xi e_{n,j}))^2=\Id_n+\xi e_{1,j},\\
u^{-1}(sr)^2u&=(\sigma r)^2 =((\Id_n+\xi e_{n,j})(\Id_n+e_{1,n}))^2=\Id_n+\xi e_{1,j}.\\
\end{array}
$$
The case $j=1$, $n\neq 1,\,n$ can be treated similarly. 
\epf

\begin{lema}\label{lema:gap}  The  unipotent classes  of type $(3)$ in $\PSL_3(2^{2m})$ are of type D for every $m\ge1$.
\end{lema}
\pf  Let $\zeta\in \F_4^\times - \F_2$, $ r= \left(\begin{smallmatrix}
1&1 &0\\
0&1 &1 \\
0&0 &1 
\end{smallmatrix}\right)$, $s= \left(\begin{smallmatrix}
\zeta^2 & 0 & \zeta^2\\
\zeta & 1 & \zeta^2 \\
\zeta & \zeta^2 & \zeta^2 
\end{smallmatrix}\right)$. Then $(rs)^2 \neq (sr)^2$, $\oc = \oc^{\SL_3(4)}_r 
= \oc^{\SL_3(4)}_s$ and $\oc^{H}_r \neq \oc^{H}_s$,
where $H = \langle r, s\rangle < \SL_3(4)$. Indeed, $|H|=108$ 
and it can be presented as the
group generated by two elements $r,s$ satisfying the relations  
$r^{4}=s^{4} = 1$, $(rs)^3 = 1$,  $(r\trid(s^{-1}\trid (r\trid s)))s^{-1}=1$.
Thus, also $(sr)^{3}=1$. In particular, 
$sr^{-1}s\in C_{H}(r)$, $rs^{-1}r\in C_{H}(s)$ and
\begin{align*}
\oc^{H}_r  & = \{ (r^{i}s^{j})\trid r ~|~ 0\leq i \leq 3, 0\leq j\leq 2 \text{ and }
i=0 \text{ if }j=0\} \text{ and }\\
\oc^{H}_s  & = \{(s^{i}r^{j})\trid s ~|~ 0\leq i \leq 3, 0\leq j\leq 2 \text{ and }
i=0 \text{ if }j=0 \}, 
\end{align*}
with $|\oc^{H}_r|=9=|\oc^{H}_s|$. A direct computation shows that 
$s\not\in\Oc_r^H$.
The Lemma follows, as $\SL_3(4) < \SL_3(2^{2m})$.
\epf

\subsection{Unipotent conjugacy classes of type F}
\label{subsec:slnq-type4}
Here assume that $q$ is even and investigate when a unipotent class is of type F;
recall that not all classes are of type D, see Proposition \ref{prop:2111-notD}.

\begin{lema}\label{lem:rack-type4-SLn}
Let $u \in G$  of
type $(\lambda_{1},\ldots ,\lambda_{k})$ and assume that either
\begin{enumerate}\renewcommand{\theenumi}{\alph{enumi}}\renewcommand{\labelenumi}{(\theenumi)}
	
	\item\label{item:typeF-even-cinco}  $\lambda_{1} \geq 5$.

	\item\label{item:typeF-even-tres-dos}  $\lambda_{1} = 3$, $\lambda_{2} = 2$.

	\item\label{item:typeF-tres} $\lambda_{1} = 3$ and $q\geq 8$, or

	\item\label{item:typeF-dostresunos} $\lambda_{1} = 2$ and $\lambda_{j} = 1$ for at least 3 different $j$.
\end{enumerate}

Then the conjugacy class $\oc_{u}$ is of type F.
\end{lema}

\pf By Lemma \ref{lem:blockD}, it is enough to look at some specific unipotent classes; when these are regular
 we may assume that $\Oc = \Oc_{r_{\uno}}$ by Remark \ref{obs:conj-class-unip-Fq}.

\begin{step}\label{lem:type4-regular5}
If $n>4$, then a regular unipotent class $\Oc$ in $G$ is of type F.
\end{step}

Let $\ba = (a_{1},a_{2}, a_{3})$, $\bu = (u_{1},u_{2}, u_{3})  \in\kc^{3}$ and set
\begin{align*}
x_{\ba}(\bu) &= \left( \begin{smallmatrix}
         1 & 1 & 0 & 0 & \ldots & a_1 & u_1  & u_3 \\
         0 & 1 & 1 & 0 & \ldots  & 0   & a_2 & u_2 \\
         0 & 0 & 1 & 1 & \ldots  & 0   & 0 & a_3 \\
       \vdots & & \vdots & \vdots  & \vdots & & \vdots & \vdots \\
         0 &\ldots & \ldots & &\ldots & 0 & 1 & 1 \\
         0 &\ldots & \ldots &\ldots & \ldots & \ldots & 0 & 1
              \end{smallmatrix}\right).
\end{align*}

Let $X_{\ba} = \left\{x_{\ba}(\bu): \bu = (u_{1},u_{2}, u_{3}) \in\kc^{3} \right\}\subset\oc$.
 Then $x_{\ba}(\bu)x_{\bb}(\bv) = x_{\bb}(\bv +\bw) x_{\ba}(\bu)$,
where $\bw = $
\begin{align*} & (a_{1} + a_{2} + b_{1} + b_{2}, a_{2} + a_{3} + b_{2} + b_{3}, a_{1} + a_{2} + b_{1} + b_{2} + u_{1} + u_{2} + v_{1} + v_{2}).
\end{align*}
Thus  $X_{\ba}\trid X_{\bb} = X_{\bb}$, for every $\ba, \bb \in \kc^{3}$. 
Let $$A = \{(1,1,1), (1,1,0), (1,0,1), (0,1,1)\}.$$ If $\ba \neq \bb \in A$, then 
$x_{\ba}(\bu)\trid x_{\bb}(\bv) \neq x_{\bb}(\bv)$ for any $\bu$, $\bv$, and $\oc$ is of type F.

\begin{step}\label{lem:type4-tres-dos}
If  $\oc$ is unipotent of type $(3,2)$, then $\oc$ is of type F.
\end{step}

Let $\ba = (a_{1},a_{2}, a_{3})$, $\bu = (u_{1},u_{2}, u_{3})  \in\kc^{3}$ and set
$x_{\ba}(\bu) = \left( \begin{smallmatrix}
         1 & 1  & a_1 & u_1  & u_3 \\
         0 & 1 & 1 &  a_2 & u_2 \\
         0 & 0 & 1 & 0 & a_3 \\
         0 &0 & 0 & 1 & 1 \\
         0 & 0 & 0 & 0 & 1
              \end{smallmatrix}\right)$.
Let $X_{\ba} = \left\{x_{\ba}(\bu): \bu = (u_{1},u_{3}, u_{3}) \in\kc^{3} \right\}$.  It can be shown that 
$X_{\ba} \subset \oc$ if and only if $\ba\in I = \{(a_{1},a_{2}, a_{3})\in\kc^{3}: a_{2} = a_{3}\}$. Let $\ba, \bb\in I$.
Now
\begin{align}\label{eq:def-Xabis}
&x_{\ba}(\bu)x_{\bb}(\bv) = x_{\bb}(\bv +\bw) x_{\ba}(\bu), \end{align}
where $\bw =  (a_{2} + b_{2}, 0,   a_{2}  + b_{2} + a_{1}b_{2} + a_{2}b_{1} + u_{1} + u_{2} + v_{1} + v_{2})$.
By \eqref{eq:def-Xabis},  $X_{\ba}\trid X_{\bb} = X_{\bb}$, for every $\ba, \bb \in I$. 
Let 
\begin{align*}
A &= \{\ba_1 = (1,0,0), \ba_2 = (1,1,1),  \ba_3 = (0,1,1), \ba_4 = (0,0,0)\} \subset I; \\
 r_1 & = x_{\ba_1}(1,0,0), \quad r_j = x_{\ba_j}(0,0,0), \quad 2\le j \le 4.
\end{align*}

Then $r_j\in R_j := X_{\ba_j}$ and $r_i\trid r_j \neq r_j$ for $i\neq j \in\I_4$, so 
$\oc$ is of type F.

\begin{step}\label{lem:type4-sln}
If $n = 3$ and $q \geq 8$, then a regular unipotent class $\Oc$ is of type F.
\end{step}

Let $A= \{\ba := (a^{2}, a^{-1}) \in (\F^{\times}_{q})^{2}: a\in \F_{q}^{\times}\}$. Then
$(R_{\ba})_{\ba \in A}$ is a family of mutually disjoint subracks of $\Oc$, 
by Lemma \ref{lem:conj-class-unip-Fq}. 
Now $\theta^{1}_{\ba,\bb} =0 \Leftrightarrow a^3 =b^3 \Leftrightarrow a =b$. Hence $r_{\ba}\trid r_{\bb} \neq r_{\bb}$ 
for $a\neq b$, by \eqref{eq:rack-ra-trid}.
As $|\F_{q}^{\times}| \geq 4$, $\Oc$ is of type F.

\begin{step}
If $u \in G$ is unipotent of type $(2,1,1,1)$,  then $\oc = \oc_{u}$ is of type F.
\end{step}

We may assume
that $u = \left(\begin{smallmatrix}
r_{1} & 0 \\
0 & \id_3
\end{smallmatrix}\right)$.
Let $(\be_{j})_{j \in \I_4}$ be the canonical basis of $\kc^{4}$ and 
$R_j = R_{\be_{j}} \cap \oc$; then $R_{j} \trid R_{k} \subseteq R_{k}$
for $k,j \in \I_4$. Let $r_1 =r_{\be_{1}}$,
$r_2 =r_{\be_{2}}$,
\begin{align*}
r_{3}& =\left(\begin{smallmatrix}
1&0&0&0&0\\
0&1&0&1&0\\
0&0&1&1&0\\
0&0&0&1&0\\
0&0&0&0&1\\
\end{smallmatrix}\right) = \left(\begin{smallmatrix}
\id_2& e_{2,1}\\
0&\id_3 \end{smallmatrix}\right)r_{\be_{3}}\left(\begin{smallmatrix}
\id_2& e_{2,1}\\
0&\id_3 \end{smallmatrix}\right), \\
r_{4}& =\left(\begin{smallmatrix}
1&0&0&0&0\\
0&1&0&0&1\\
0&0&1&0&1\\
0&0&0&1&1\\
0&0&0&0&1\\
\end{smallmatrix}\right) = \left(\begin{smallmatrix}
\id_3& e_{2,1} + e_{3,1}\\
0&\id_2 \end{smallmatrix}\right)r_{\be_{4}}\left(\begin{smallmatrix}
\id_3& e_{2,1} + e_{3,1}\\
0&\id_2 \end{smallmatrix}\right).
\end{align*}
Then $r_{j} \in R_{j}$ and
$r_{j}\trid r_{k} \neq r_{k}$, $j\neq k \in \I_4$. Thus
$\Oc$ is of type F.
\epf

By Proposition  \ref{prop:2111-notD}, the classes of type $(2,1)$ in $\SL_3(q)$, $q$ even, are
not of type D. Now we show that they are not of type F, hence are chtulhu.

\begin{prop}\label{prop:21-notF} The unipotent classes of type $(2,1)$ in $G = \SL_3(q)$  for $q$ even are not of type F. \end{prop}

\pf Let $(r_a)_{a\in\I_4}$ in $G$ such that
$\Oc_{r_a}^{\langle r_a: a\in \I_4\rangle}\neq \Oc_{r_b}^{\langle r_a: a\in \I_4\rangle}$
and $r_a\neq r_b$, for all $a\neq b$ in $\I_4$. 
Without loss of generality we may assume $r_1= \Id_3 + e_{1,3}$. Then, arguing as in the proof of  Proposition \ref{prop:2111-notD}
we have:
\begin{itemize} 
\item $r_a=u_an_at_a\triangleright r_1$, for $a\in\{2,3,4\}$, where $u_a\in\U^F$, $t_a\in\T^F$ and $n_a$ monomial in $\SL_3(2)$;
\item $\sigma_a:=n_at_a\triangleright r=\id_3+\xi_a e_{i_a,j_a}$ for some $\xi_a\in\F_q^\times$ and
$(i_a,\,j_a)\in\{(2,1),\,(3,\,2)\}$. 
\end{itemize}
Thus,  there are  $a\neq b$ in $\{2,3,4\}$ such that
$(i_a,\,j_a)=(i_b,\,j_b)$.  We claim that if $(i_a,\,j_a)=(i_b,\,j_b)=(2,\,1)$
then $|r_ar_b|$ is either $2$ or odd. Hence by Remark  \ref{obs:rack-dihedral}, 
either $r_ar_b = r_br_a$ or 
$\Oc_{r_a}^{\langle r_a, r_b\rangle}= \Oc_{r_b}^{\langle r_a, r_b\rangle}$, 
a contradiction to our assumption.
Since matrices in $\id_3+\F_q e_{2,3}$ commute with $\sigma_a$ and $\sigma_b$, 
there is no loss of generality in taking 
$u_a,\,u_b\in \id_3+\F_q e_{1,2}+\F_q e_{1,3}$. Further, 
$$|r_ar_b|=|u_a\sigma_au_a^{-1}u_b\sigma_bu_b^{-1}|=
|\sigma_a((u_a^{-1}u_b)\triangleright\sigma_b)|$$ so to prove the claim we may take
$r_a=\sigma_a$, $r_b=(u_a^{-1}u_b)\triangleright\sigma_b$. 
Then, for $u_a^{-1}u_b = \id_3 + xe_{1,2}+ye_{1,3}$ we have
\begin{align*}
&& r_b &= \left(\begin{smallmatrix}1+\xi_bx&\xi_bx^2&\xi_bxy\\
\xi_b&1+\xi_b x&\xi_by\\
0&0&1\end{smallmatrix}\right),& r_ar_b &= \left(\begin{smallmatrix}A&{\mathbf c}\\ 
0&1\end{smallmatrix}\right),\\ &\text{where}&
A &= \left(\begin{smallmatrix}1+\xi_bx&\xi_bx^2\\
\xi_a+\xi_a\xi_bx+\xi_b&\xi_a\xi_b x^2+1+\xi_bx\end{smallmatrix}\right),& {\mathbf c} &=\left(\begin{smallmatrix}\xi_bxy\\
\xi_a\xi_bxy+\xi_by\end{smallmatrix}\right).
\end{align*}

Now, $(r_ar_b)^k=\left(\begin{smallmatrix}A^k&(A^{k-1}+\cdots+\id_2){\mathbf c}\\
0&1\end{smallmatrix}\right)$. 
Besides, $A\in\SL_2(q)$ so it is either semisimple or unipotent, the latter occurring if and only if $Tr(A)=0$,  if and only if $x=0$. 
In this case, $r_ar_b = r_br_a$. Otherwise $A$ is semisimple, hence 
$|A|=h$ is odd and $A^{h-1}+\cdots+\id_2 =0$, so $|r_ar_b|=h$; the claim follows. 
\epf

\subsection{Collapsing unipotent classes in $\Gb = \PSL_{n}(q)$}\label{subsec:summarize-collapse}
We summarize the results in \S \ref{subsec:slnq-typeD} and \ref{subsec:slnq-type4}  
showing the unipotent classes in $\Gb$  that collapse in  \ref{tab:unipotent-chevalley-collapse-typeA}.
%We apply systematically Lemma \ref{lem:blockD}.
Recall that we assume $q\neq 2$, $3$, $4$, $5$, $9$, when $n=2$.
The information in \ref{tab:unipotent-chevalley-collapse-typeA} is minimal; 
many orbits collapse by different reasons, but we omit to discuss this in detail.

\begin{align}\label{tab:unipotent-chevalley-collapse-typeA}\stepcounter{tabla}\tag*{Table \thetabla}
\begin{tabular}{|c|c|c|c|}
\hline $n$   & $q$ & type $(\lambda_1, \dots, \lambda_k)$ & Criterium  \\
\hline
$2$   & odd square $> 9$      & $(2)$ &  D,    \ref{lem:sl2} \\
 \hline\hline
$ > 2$  & odd   & $\lambda_{1} \ge 3$ &  D,  \ref{cor:rack-typeD-SLn} \\
  \cline{3-4}
   &    &  $(2, 2, \dots)$ &  D,  \ref{lem:2-1,2-2} \eqref{item:2-2} \\
  \cline{3-4}
   &    &  $(2, 1 \dots)$ &  D,  \ref{lem:2-1,2-2} \eqref{item:2-1} \\
 \hline\hline
 &  even & $\lambda_1 \ge 5$ & F,  \ref{lem:rack-type4-SLn} \eqref{item:typeF-even-cinco}
\\
 \cline{3-4}
 &   & $\lambda_1 = 4$ & D,  \ref{lema:q=2-n=4} \eqref{item:typeD-even-cuatro}
\\ \cline{3-4}
  &     & $(3, 3, \dots)$  & D, \ref{lem:difblocksallq}
\\ \cline{3-4}
  &     & $(3, 2, \dots)$  & F,   \ref{lem:rack-type4-SLn} \eqref{item:typeF-even-tres-dos}
 \\
 \cline{3-4}
 &   & $(3, 1, \dots)$ & D,   \ref{lema:q=2-n=4} \eqref{item:typeD-even-tres-uno}
\\ \cline{3-4}
 &   & $(2, 2, \dots)$ & D,   \ref{lema:q=2-n=4} \eqref{item:typeD-dos-dos}
\\ 
\cline{3-4}
 &   & $(2,1,1,1,\dots)$ &  F,  \ref{lem:rack-type4-SLn} \eqref{item:typeF-dostresunos}
\\
\hline
& even $\geq 8$ & $\lambda_1 = 3$ &   F,  \ref{lem:rack-type4-SLn} \eqref{item:typeF-tres}\\
  & 4 &   &D, \ref{lema:gap} 
\\
%\cline{4-4}\cline{2-2}
 \hline
 \end{tabular}
 \end{align}
 
\bigbreak
Now we deal with the Nichols algebras of irreducible Yetter-Drinfeld modules associated
to the remaining classes in  \ref{tab:uno}. We recall the useful
\emph{little triangle} Lemma. 
Let $G$ be a finite group. A conjugacy class $\Oc$ in $G$ contains a
\emph{little triangle} if there are different elements $(\sigma_i)_{i\in \I_3}$  such that
\begin{itemize}
\item $\sigma_1^{h}=\sigma_2\sigma_3$ for an odd integer $h$;
\item $\sigma_i\sigma_j=\sigma_j\sigma_i$, $i,j\in \I_3$;
\item there are $g_2$, $g_3\in G$ such that
$\sigma_{i}=g_{i}\sigma_{1}g_{i}^{-1}$ and $g_{3}g_{2}, g_{2}g_{3}\in C_{G}(\sigma_1)$.
\end{itemize}

\begin{lema}\label{lem:triangulitos} \cite[Lemma 2.3]{FGV2}
If $\Oc$ in $G$ contains a little triangle, then $\dim\toba(\Oc,\rho)=\infty$, for every
$\rho\in \Irr C_{G}(\sigma_1)$. \qed
\end{lema}

Clearly, if $\Oc^G_x$ of $ x\in G$ contains a little triangle and $\psi: G\to H$ is a group homomorphism,
then $\Oc^H_{\psi(x)}$ contains a little triangle. In particular,
\begin{itemize}
\item If $x\in G < H$ and $\Oc_x^G$ of $G$ contains a little triangle, then  $\Oc_x^H$  also contains a little triangle.
\item If $T$ is an (outer) automorphism and $\Oc$ in $G$ contains a little triangle, then so does $T(\Oc)$.
\end{itemize}

\begin{lema}\label{lem:sl2-YD-modules}
Let $\Oc$ be the conjugacy class of $x\in \Gb$. In the cases listed below, $\dim\toba(\Oc,\rho)=\infty$, for every
$\rho\in \Irr C_{\Gb}(x)$.

\begin{enumerate}\renewcommand{\theenumi}{\alph{enumi}}\renewcommand{\labelenumi}{(\theenumi)}
	\item \label{step:sl2-YD-dos} $\Gb = \PSL_2(q)$, $x$ of type $(2)$.
	
	\item \label{step:sl2-YD-dos-uno} $\Gb = \PSL_3(q)$ with $q$ even, $x$ of type $(2,1)$.
	
	\item \label{step:sl2-triangulito-211} $\Gb = \PSL_4(q)$  with $q$ even, $x$ of type $(2,1,1)$.
	
\end{enumerate}
\end{lema}

\pf   \eqref{step:sl2-YD-dos}: \cite[3.1]{FGV1} for $q$ even, \cite[4.1, 4.3]{FGV2} for $q$ odd.
\eqref{step:sl2-YD-dos-uno} and \eqref{step:sl2-triangulito-211}:
$\PSL_3(2) \simeq \PSL_2(7)$ contains a copy of $\ac$, so the class of involutions contains a little triangle \cite[4.3]{FGV2}.
Now the previous remarks apply.
\epf

\section{Non-semisimple classes in $\PSL_{n}(q)$}\label{sec:non-ss}
\subsection{Preliminaries} \label{subsec:nonss-prelim}
In this section we apply the results in Section \ref{sec:sln} on unipotent classes to
non-semisimple classes in $\Gb = \PSL_{n}(q)$.
Let $x\in \Gb$ and pick $\x\in \SL_{n}(q)$ such that $\pi(\x) = x$; if $\x =
\x_s\x_u$ is  the Chevalley-Jordan decomposition of $\x$, then
$x_s = \pi(\x_s)$ and $x_u = \pi(\x_u)$ form the Chevalley-Jordan
decomposition of $x$. Now $\x_u$ belongs to $\K :=
C_{\SL_{n}(q)}(\x_s)$, thus $x_u\in K := \pi(\K)$ and there are
morphisms of racks $\oc^{\K}_{\x_u}\simeq  \oc^{K}_{x_u}
\hookrightarrow \Oc_x^{\Gb}$. Hence, in many cases it will be enough to deal with
$\oc^{\K}_{\x_u}$ and to start with we describe $\K =
C_{\GL_{n}(q)}(\x_s) \cap \SL_{n}(q)$. Up to conjugation by a
matrix in $\SL_{n}(q)$, we may assume that
\begin{equation}\label{eq:u-semisimple-sln}
\x_s =
\left( \begin{smallmatrix}
         \St_{1} & 0 &    \ldots  & 0 \\
         0 & \St_{2} &  \ldots & 0 \\
       \vdots & & \ddots  & \vdots \\
                 0 &\ldots & \ldots & \St_{k}
              \end{smallmatrix}\right),
 \end{equation}
with $\St_{i}  \in \GL_{\lambda_{i}}(q)$ irreducible, that is,
its characteristic polynomial $\chi_{\St_i}$ is irreducible in $\F_q[X]$. 
Furthermore, $\prod_{i\in \I_k} \det \St_i = 1$. 
Now, if  $\sigma\in \s_k$, 
then there is $T\in \SL_{n}(q)$ such that
$T\x_sT^{-1} =\left( \begin{smallmatrix}
         \St_{\sigma(1)} & 0 &    \ldots  & 0 \\
         0 & \St_{\sigma(2)} &  \ldots & 0 \\
       \vdots & & \ddots  & \vdots \\
                 0 &\ldots & \ldots & \St_{\sigma(k)}
              \end{smallmatrix}\right)$.

\medbreak
If $\St\in \GL_{\Lambda}(q)$ is irreducible, then the subalgebra $\C_\St$ of
matrices commuting with $\St$  is a division ring by Schur Lemma;
being finite, is isomorphic to $\F_{q^{\mu}}$ for some $\mu \in \N$.
We claim that $\mu = \Lambda$. Indeed, the characteristic and minimal polynomials of $\St$ coincide
and have degree $\Lambda$, so standard arguments for finite fields imply the claim.

\begin{obs}\label{rem:regroup}
Let $\St$, $\Rt \in \GL_{\Lambda}(q)$ be semisimple and conjugate in $\GL_{\Lambda}(\kk)$. Then there exists $T\in\SL_{\Lambda}(q)$ 
such that $T \St T^{-1} = \Rt$; that is, $\St$ and $\Rt$ are conjugate under $\SL_{\Lambda}(q)$. 
\end{obs}

Indeed, $\Rt$ and $\St$ are conjugate in $\GL_{\Lambda}(q)$ by \cite[8.5]{Hu}, \cite[I.3.5]{sp-st}. Also we may assume that $\St$ is irreducible.
Let $T_0\in\GL_{\Lambda}(q)$ such that $T_0 \St T_0^{-1} = \Rt$. Since $\det:\C_{\St}^{\times} \to \F_q^{\times}$ equals the norm 
$N: \F_{q^{\Lambda}}^{\times} \to \F_q^{\times}$ which is surjective, we may pick $T_1\in \C_{\St}$ such that 
$\det T_1 = \det T_0^{-1}$. Then $T = T_0T_1$ does the job. 

\medbreak
Assume that $\St$ is irreducible but not in $\F_q$; then $\chi_\St(\St^q) =\left(\chi_\St(\St)\right)^q = 0$, so $\St$ and $\St^q$ are conjugate 
under $\SL_{\Lambda}(q)$, but $\St\neq \St^q$.
Indeed, it can be shown that $\St^{q^i}$, $i\in \I_{\Lambda}$, are all the roots of $\chi_\St$ in $\F_q[\St] \simeq \F_{q^\Lambda}$.

\begin{obs}\label{rem:quasi-real-psl}
Let $\pi: \GL_{\Lambda}(q) \to \PGL_{\Lambda}(q)$ and let $\St \in \GL_{\Lambda}(q)$ irreducible with $\Lambda > 1$;
hence $\St \neq \St^q$.  
Then $\pi(\St) = \pi(\St^q)$ if and only if $\chi_\St$ belongs to 
\begin{align}\label{eq:irr-fact}
\irrfact(q) = \{F \in \F_q[X] &\text{ irreducible}: F \vert X^{q-1} - c, \\
\notag &\text{ for some } c\in \F_q^{\times},\, c\neq 1, \, c^{\deg F} =1\}.
\end{align}

\end{obs}

\subsection{Centralizers} \label{subsec:nonss-centralizer}
By the previous considerations, we may regroup the blocks so that there exist integers $h_1, \dots, h_\ell$
such that $\St_{i}$ and $\St_{j}$ are conjugate under $\SL_{\lambda_{i}}(q)$   
if and only if there exists a (unique) $t\in \I_{\ell}$ such that
$i, j \in \J_t$, where
\begin{align}\label{eq:blocks-ss}
\J_t&=\{i\in \N: h_1+ \dots + h_{t-1} + 1 \leq i \leq h_1+ \dots +h_t\}.
\end{align}
So, we set $\Lambda_t = \lambda_i$, if $i \in \J_t$,
$t\in \I_{\ell}$. In other words,
$h_1$ is the number of blocks $\St_i$ that are isomorphic to $\St_1$, all
of size $\Lambda_1$; $h_2$ is the number of blocks $\St_i$ that are isomorphic to $\St_{h_1 + 1}$, all
of size $\Lambda_2$, and so on.

\begin{prop}\label{prop:carter-sln} \cite{carter-classical}
$C_{\GL_{n}(q)}(\x_s)\simeq \GL_{h_1}(q^{\Lambda_1}) \times \dots
\times \GL_{h_\ell}(q^{\Lambda_\ell})$.
\end{prop}

\pf Let $\St \in \GL_{N}(q)$, $\Rt\in \GL_{P}(q)$ be irreducible.
Let $Z = \begin{pmatrix} A & B\\ C & D \end{pmatrix} \in M_{N+
P}(q)$, where $A$ is of size $N \times N$. Then $Z$ commutes with
$\begin{pmatrix} \St & 0\\ 0 & \Rt \end{pmatrix}$ iff
\begin{align*}
A\St &= \St  A, & B\Rt&= \St B,& C\St&= \Rt C, & D\Rt &= \Rt D.
\end{align*}
So that $A\in \C_{\St} \simeq \F_{q^{\Lambda}}$, $D\in \C_{\Rt} \simeq
\F_{q^{\nu}}$. If $\St$ and $\Rt$ are not conjugated, then $B=0$, $C=0$ by Schur Lemma.
Otherwise, $N = P$; we may assume $\St = \Rt$, hence $A, B, C, D \in \C_{\St} \simeq
\F_{q^{\Lambda}}$. The claim follows from this.
For, assume that
$\x_s$ is of the form \eqref{eq:u-semisimple-sln}.
Let $Z = (Z_{ij}) \in \GL_n(q)$,
where $Z_{ij}\in \F_q^{\lambda_i \times \lambda_j}$, $i,j \in \I_k$. Then
$Z\in C_{\GL_{n}(q)}(\x_s)$ iff
$Z_{ij} =0$ unless $i,j\in \J_t$ for some $t$,
in which case $Z_{ij} \in \F_{q^{\Lambda_t}}$.
Thus every $Z\in C_{\GL_{n}(q)}(\x_s)$ is a matrix of blocks
\begin{align}\label{eq:blocks-isom} Z =
\left( \begin{smallmatrix}
         W_{1} & 0 &    \ldots  & 0 \\
         0 & W_{2} &  \ldots & 0 \\
       \vdots & & \ddots  & \vdots \\
                 0 &\ldots & \ldots & W_{\ell}
              \end{smallmatrix}\right)
\end{align}
where in turn $W_t$ is a matrix of $h_t^2$ blocks,
each of size $\Lambda_t\times \Lambda_t$ and
belonging to $\C_{\St_i}\simeq\F_{q^{\Lambda_t}}$,
if $i\in \J_t$, $t\in \I_\ell$.
Thus $W_t$ can be thought of as a matrix
$\widetilde{W}_t\in M_{h_t}(q^{\Lambda_t})$,
and the map $\psi_t:W_t\mapsto \widetilde{W}_t$
is an isomorphism of monoids.
Also, $\det Z \neq 0$ iff $\det W_t \neq 0$ in
$\GL_{h_t\Lambda_t}(q)$ for all $t\in\I_\ell$,
iff $\det \widetilde{W}_t \neq 0$ in $\GL_{h_t}(q^{\Lambda_t})$
for all $t\in\I_\ell$.
Thus $\psi_t$ gives rise to an isomorphism from
the group $G_t$ of matrices
\eqref{eq:blocks-isom} with all $W_r = \id$,
except for $r=t$, to $\GL_{h_t}(q^{\Lambda_t})$.
\epf

Let $\Psi: C_{\GL_{n}(q)}(\x_s)\to \GL_{h_1}(q^{\Lambda_1}) \times \dots
\times \GL_{h_\ell}(q^{\Lambda_\ell})$ be the isomorphism given by Proposition \ref{prop:carter-sln}. Then
\begin{align}\label{eq:structure-K}
\K \simeq \{X \in \GL_{h_1}(q^{\Lambda_1}) \times \dots
\times \GL_{h_\ell}(q^{\Lambda_\ell}): \det\circ \Psi^{-1} (X) = 1\}.
\end{align}
In particular, if $\SL_{h_t}(q^{\Lambda_t}) \neq \SL_{2}(2)$, $\SL_{2}(3)$, then it is perfect, hence
\begin{align}\label{eq:K-contiene}
\SL_{h_t}(q^{\Lambda_t}) &\hookrightarrow\K,& 1\le t&\le \ell.
\end{align}
If $\SL_{h_t}(q^{\Lambda_t}) = \SL_{2}(2)$ or $\SL_{2}(3)$, then 
\eqref{eq:K-contiene} also holds, being $\Lambda_t = 1$.

\subsection{End of the proof of Theorem \ref{th:unipotent-slnq-collapse}} \label{subsec:nonss-mainthm}

\begin{prop}\label{prop:nonss-typeDPSLn}
Let $x\in\Gb$ be neither semisimple nor unipotent. Then $\oc^{\Gb}_{x}$ collapses.
%Let $\x\in \SL_{n}(q)$ with Chevalley-Jordan
%decomposition $\x = \x_s\x_u$,
%with $\x_s$ not central and $\x_u\neq e$.
%Then 
%$\oc^{\K}_{\x_u}$ collapses.
%In consequence, if  
%If $x = \pi(\x) \in \Gb$, then $\oc^{\Gb}_{x}$ collapses.
\end{prop}

\pf Let $\x\in \SL_{n}(q)$ with $x = \pi(\x)$, and let $\x = \x_s\x_u$ be its
Chevalley-Jordan decomposition. By our assumption $\x_s$ is not central and $\x_u\neq e$.

We assume that $\x_s$ is in the form \eqref{eq:u-semisimple-sln}; then there are natural numbers $h_1, \dots, h_\ell$,
$\Lambda_1, \dots, \Lambda_\ell$ such that the structure of $\K$
is given by \eqref{eq:structure-K}. Then $\x_u = (u_1, \dots, u_\ell)$
with $u_t\in \GL_{h_t}(q^{\Lambda_t})$ unipotent, $ t\in \I_{\ell}$.
For simplicity, we write also $\x_s = (S_1, \dots, S_\ell)$.
Up to a further reordering, there exists $M\in \I_{\ell}$ such that $u_{t}\neq \id$ iff $t \leq M$,
and $h_1 \geq \dots \geq h_M > 1$; since $\x_u\neq e$, $M > 0$.
Recall that $\oc^{\SL_{h_t}(q^{\Lambda_t})}_{u_t}$ is a subrack of
$\oc^{\K}_{\x_u}$ for all $t$ by \eqref{eq:K-contiene}.

\smallbreak
By the unipotent part of Theorem \ref{th:unipotent-slnq-collapse}, we may assume
that  $h_t \leq 4$ and  $u_t$ appears in \ref{tab:uno} or it is of
type $(2)$ and $q$ is in $\{2,3,4,5,9\}$, for all $t\in \I_M$.
Let $X$ be a unipotent orbit either of type $(3)$ with $q^{\Lambda}= 2$; or else of type $(2)$ with $q^{\Lambda}$ even or 9 or odd not a square;
or else of type $(2,1)$ or $(2,1, 1)$ with $q^{\Lambda}$ even.
By inspection, we see that

\begin{enumerate}\renewcommand{\theenumi}{\alph{enumi}}\renewcommand{\labelenumi}{(\theenumi)}
  \item\label{eq:a} There exist $x_1, x_2 \in X$ such that $(x_1x_2)^2  \neq (x_2x_1)^2$.
  \item\label{eq:b} There exist $y_1, y_2 \in X$ such that $y_1y_2 = y_2 y_1$, except when $X$ is of type $(2)$ with $q^{\Lambda} = 2$ or $3$.
\end{enumerate}

\begin{stepnss}\label{case:nonss-uno}
$M = 1 = \ell$. Then $\Oc_x^{\Gb}$ is of type D.
\end{stepnss}

In this case $x_u=u_1$. In addition, $\Lambda = \Lambda_{1} > 1$ since $\x_{s}$ is not central; so $q^\Lambda\neq2$. Hence type $(3)$ and $(2)$ with $q^{\Lambda}=2$ are excluded.
Let $\St = \St_1$. Assume  that   $\chi_{\St}  \notin \irrfact(q)$.
By Remarks \ref{rem:regroup} and \ref{rem:quasi-real-psl},
\begin{align*}
\left(\oc^{\K}_{\x_u}\right)^{(2)}\, \overset{\ref{lema:prod-racks-D}}= \, \oc^{\K}_{\x_u} \coprod \oc^{\K}_{\x_u} \simeq  \pi(\St)\oc^{K}_{x_u} \coprod \pi(\St^q)\oc^{K}_{x_u}
\hookrightarrow \Oc_x^{\Gb}.
\end{align*}
By \eqref{eq:a}, $\Oc_x^{\Gb}$ is of type D.
Now, if   $\chi_{\St} \in \irrfact(q)$, then $\St$ is conjugated to $\St^q = c\St$ for some $c\in \F_q^{\times} - 1$. Pick $Y\in \SL_\Lambda(q)$ such that
$Y\St Y^{-1} = c\St$. If $\x_u$ is of type $(2)$, then take 
\begin{align*}
r &= \left(\begin{smallmatrix}
\St & \St  \\
0 & \St
 \end{smallmatrix}  \right), &s & = 
\left(\begin{smallmatrix}
\id & 0  \\
0 & Y
 \end{smallmatrix} \right) \trid r = \left(\begin{smallmatrix}
\St & \St Y^{-1}  \\
0 & c\St
 \end{smallmatrix} \right); \\
R_1 &= \left\{\left(\begin{smallmatrix}
\St & * \\
 0 &\St \end{smallmatrix}\right) \in \F_q^{2\Lambda\times 2\Lambda} \right\},& R &= \pi(R_1) \cap \oc \ni \pi(r); \\
\St_1 &= \left\{\left(\begin{smallmatrix}
\St & * \\
 0 & c\St \end{smallmatrix}\right) \in \F_q^{2\Lambda\times 2\Lambda} \right\},& \St &= \pi(\St_1) \cap \oc \ni \pi(s).
\end{align*}
Then $r$ and $s$ are conjugated in $\SL_{2\Lambda}(q^\Lambda)$ and $R \coprod S \hookrightarrow \Oc_x^{\Gb}$ is decomposable. Also
$(rs)^2 = (sr)^2$ means that
\begin{align*}
&\left(\begin{smallmatrix}
\St^4\, & (c+c^2)\St^4 + \St^4Y^{-1} + c\St^2Y^{-1}\St^2   \\
0 & c^2\St^4
 \end{smallmatrix}  \right) 
 =  \left(\begin{smallmatrix}
\St^4\, & (c + 1)\St^4 + \St^3Y^{-1}\St + c\St Y^{-1}\St^3   \\
0 & c^2\St^4
 \end{smallmatrix}  \right) 
 \\ &\iff \St^4Y^{-1} + c^2\St^4 + c\St^2Y^{-1}\St^2 =\St^3Y^{-1}\St + \St^4 + c\St Y^{-1}\St^3
  \\ &\iff c^2(c^2-1)\id= (1- c^2)Y \\ & \iff c^2 =1.
\end{align*}
where we have used $\St Y^{-1}=cY^{-1}\St$ and that $Y$ is not a scalar matrix.
Thus, if $q$ is even, then $c=1$, a contradiction; and if $q$ is odd and $c\neq 1$, then $\ord c = 2$, hence $\Lambda$ is even and $q^\Lambda$ is a square.
Hence $\Oc_x^{\Gb}$ is of type D, except when $q^\Lambda =9$. If $q^\Lambda =9$, then $q=3$ and $\Lambda =2$. 
Let $\St = \left(\begin{smallmatrix}
0 & 1   \\
2 & 0
 \end{smallmatrix}  \right)$, $\Rt = \left(\begin{smallmatrix}
1 & 1   \\
1 & 2
 \end{smallmatrix}  \right)\in \SL_2(3)$; they are conjugated in $\SL_2(3)$
and $\St\Rt = -\Rt\St$, so that $\pi(\St)\pi(\Rt) = \pi(\Rt)\pi(\St)$. 
It is enough to deal with $\Oc = \Oc_x^{\Gb}$ where 
\begin{align*}
\x &= r = \left(\begin{smallmatrix}
\St & \St  \\
0 & \St
 \end{smallmatrix}  \right),& 
 s & = \left(\begin{smallmatrix}
0 & \id_2  \\
2\id_2 & 0
 \end{smallmatrix} \right) \trid \left(\begin{smallmatrix}
\Rt & \Rt  \\
0 & \Rt
 \end{smallmatrix}  \right) = \left(\begin{smallmatrix}
\Rt & 0  \\
2\Rt & \Rt
 \end{smallmatrix} \right); \\
R_1 &= \left\{\left(\begin{smallmatrix}
a\St & b\St \\
d\St & c\St \end{smallmatrix}\right) \in \SL_4(9): a,b,c,d\in \F_3 \right\},& R &= \pi(R_1) \cap \oc \ni \pi(r); \\
S_1 &= \left\{\left(\begin{smallmatrix}
a\Rt & b\Rt \\
d\Rt & c\Rt \end{smallmatrix}\right) \in \SL_4(9): a,b,c,d\in \F_3 \right\},& S &= \pi(S_1) \cap \oc \ni \pi(s).
\end{align*}
Then $r$ and $s$ are conjugated in $\SL_4(9)$, $(rs)^2 \neq (sr)^2$, $R \coprod S \hookrightarrow \Oc_x^{\Gb}$ 
and $\Oc$ is of type D. 
The types $(2, 1)$ and $(2,1,1)$ are treated as above, with
\begin{align*}
r &= \left(\begin{smallmatrix}
\St & \St &0 \\
0 & \St &0 \\
0 &0  & \St 
 \end{smallmatrix}  \right),& &\text{respectively } r = \left(\begin{smallmatrix}
\St & \St &0&0 \\
0 & \St &0 &0 \\
0 &0  & \St &0 \\
0 &0  &0 & \St  \end{smallmatrix}  \right),
\\
s & = 
\left(\begin{smallmatrix}
\id & 0 & 0 \\
0 & Y &0  \\
0& 0&  \id
 \end{smallmatrix} \right) \trid r = \left(\begin{smallmatrix}
\St & \St Y^{-1} &0 \\
0 & c\St &0 \\
0 &0  & \St 
 \end{smallmatrix} \right),& &\text{respectively } 
s  = 
 \left(\begin{smallmatrix}
\St & \St Y^{-1} &0 &0 \\
0 & c\St &0 &0 \\
0 &0  & \St &0 \\
0 &0  &0 & \St 
 \end{smallmatrix} \right).
\end{align*}

\begin{stepnss}\label{case:nonss-odd-tres}
$M  = 1 < \ell$. Then $\Oc_x^{\Gb}$ is of type D.
\end{stepnss}

In this case $x_u=(u_1,1,...,1)$. Assume that $\Lambda_i > 1$ for some $i$. Then  $\y_s = (\St_1, \dots, \St_i^q, \dots, \St_\ell)$
(all $\St_h$ equal to $\St_i$ raised to the $q$), is conjugated to $\x_s$; clearly $\pi(\x_s) \neq \pi(\y_s)$. 
By Remarks \ref{rem:regroup} and \ref{rem:quasi-real-psl}
\begin{align*}
\left(\oc^{\K}_{\x_u}\right)^{(2)}\, \overset{\ref{lema:prod-racks-D}}= \, \oc^{\K}_{\x_u} \coprod \oc^{\K}_{\x_u} \simeq  \pi(\x_s)\oc^{K}_{x_u} \coprod \pi(\y_s)\oc^{K}_{x_u}
\hookrightarrow \Oc_x^{\Gb}.
\end{align*}
By \eqref{eq:a}, $\Oc_x^{\Gb}$ is of type D. Assume then that $\Lambda_i = 1$ for all $i \in \I_\ell$.
Since $\ell>1$ the case $u_1$ of type $(3)$ with $q = 2$ is excluded, so $u_1$ is of type $(2)$, $(2, 1)$ or $(2,1,1)$. 
We consider first the case when $\ell = 2$ and $u_1$ is of type $(2)$. Let 
\begin{align*}
r &= \left(\begin{smallmatrix}
\St_1 & \St_1 & 0  \\
0 & \St_1 & 0 \\
0 & 0 & \St_3
 \end{smallmatrix}  \right), &s & = 
 \left(\begin{smallmatrix}
\St_3 & 0 & 0  \\
0 & \St_1 & \St_1 \\
0 & 0  & \St_1
 \end{smallmatrix}  \right); \\
R_1 &= \left\{\left(\begin{smallmatrix} \St_1 &* & *  \\
0 & \St_1 & * \\
0 & 0 &\St_3 \end{smallmatrix}\right) \in \F_q^{3\times 3} \right\},& R &= \pi(R_1) \cap \oc \ni \pi(r); \\
\St_1 &= \left\{\left(\begin{smallmatrix} \St_3 &* & *  \\
0 & \St_1 & * \\
0 & 0 &\St_1 \end{smallmatrix}\right) \in \F_q^{3\times 3} \right\},& \St &= \pi(\St_1) \cap \oc \ni \pi(s).
\end{align*}
Then $r$ and $s$ are conjugated in $\SL_3(q)$, $R \coprod S \hookrightarrow \Oc_x^{\Gb}$ is decomposable,
\begin{align*}
(rs)^2 &=   \left(\begin{smallmatrix}
\St_1^2\St_3^2\, & \St_1^3(\St_1 + \St_3) \,& \St_1^3(\St_1 + 2\St_3)  \\
0 & \St_1^4 & \St_1^3(\St_1 + \St_3) \\
0 & 0  & \St_1^2\St_3^2
 \end{smallmatrix}  \right) 
 \\ 
(sr)^2 &=  \left(\begin{smallmatrix}
\St_1^2\St_3^2\, & \St_1^2\St_3(\St_1 + \St_3) \, & \St_1^2\St_3^2  \\
0 & \St_1^4 & \St_1^2\St_3(\St_1 + \St_3) \\
0 & 0  & \St_1^2\St_3^2
 \end{smallmatrix}  \right).
\end{align*}
Hence $(\pi(r)\pi(s))^2 \neq (\pi(s)\pi(r))^2$ and thus $\Oc_x^{\Gb}$ is of type D.
The other cases are dealt with in a similar way.

\begin{stepnss}\label{case:nonss-odd-uno}
$M >1$, and $q^{\Lambda_{t}}\neq 3$ for some $t\in \I_M$. Then $\Oc_x^{\Gb}$ is of type D.
\end{stepnss}

Assume $q$ is odd. Since $M>1$ there is $k\in \I_M - \{t\}$ such that $u_k\neq1$. We set $X =  \oc^{\SL_{h_k}(q^{\Lambda_{k}})}_{u_{k}}$, $Y = \oc^{\SL_{h_t}(q^{\Lambda_{t}})}_{u_{t}}$.
By  Lemma \ref{lema:prod-racks-D}, \eqref{eq:a} and \eqref{eq:b} 
%$\oc^{\SL_{2}(q^{\Lambda_{k}})\times \SL_{2}(q^{\Lambda_{t}})}_{(u_{k},u_{t})}$
$X\times Y$, and $\oc^{\K}_{\x_u}$, are of type D.

If $q$ is even, then the same argument applies except when $u_t$ is of type $(2)$ with $q^{\Lambda_{t}} = 2$ for all $t\in \I_M$. But here $\St_t\in\F_2^{\times}$, i. e., $\St_t = 1$ so $M =1$, a contradiction. 

\begin{stepnss}\label{case:nonss-odd-dos}
$M >1$, and $q^{\Lambda_{t}} = 3$ for all $t\in \I_M$. Then $\Oc_x^{\Gb}$ is of type D.
\end{stepnss}

According to our reduction we need only to consider the case in which $h_t=2$ and $u_t$ is of type $(2)$ for every $t\in \I_M$. 
Then $\St_t\in \C_{S_t}\simeq\F_{q^{\Lambda_t}} = \F_3$, so $\St_1 = \St_2 = 1$, $\St_3 = \St_4 = 2$ and $M =2$. 
The rack of unipotent conjugacy classes in  $\SL_2(3)$ is the union of two conjugacy
classes $\oc_1 \ni r_1$ and $\oc_2$, both isomorphic to the tetrahedral rack. 
If $M < \ell$, then 
$\K \supseteq \{(g_1, g_2, g_3) \in \GL_2(3) \times \GL_2(3) \times \F_3^{\times}: \det g_1\det g_2 = g_3^{-1}\}$; thus
$\oc^{\K}_{(r_1,r_1)} =  \coprod_{i, j \in \I_2} \oc_i \times \oc_j $ is of type D, being $\oc_1 \times \oc_1 \coprod \oc_1 \times \oc_2$ 
of type D.

If $M = \ell = 2$, then $\oc^{\K}_{(r_1,r_1)}$ is cthulhu, so we consider $\oc = \oc^{\Gb}_{x}$.
There are two different conjugacy classes with the same $x_s$, namely those with representative
$r = \left(\begin{smallmatrix}
1 & 1 &0 & 0  \\
0 & 1 &0 & 0 \\
0 & 0 &2 & 2 \\
0 & 0 & 0 & 2
 \end{smallmatrix}  \right)$, respectively $\left(\begin{smallmatrix}
1 & 2 &0 & 0  \\
0 & 1 &0 & 0 \\
0 & 0 &2 & 2 \\
0 & 0 & 0 & 2
 \end{smallmatrix}  \right)$, but they are isomorphic as racks being conjugated in $\PGL_3(4)$.
Now let $s = \left(\begin{smallmatrix}
1 & 0 &0 & 0  \\
0 & 0 &1 & 0 \\
0 & 0 &0 & 1 \\
0 & 1 & 0 & 0
 \end{smallmatrix}  \right) \trid r = 
 \left(\begin{smallmatrix}
1 & 0 &0 & 1  \\
0 & 2 &2 & 0 \\
0 & 0 &2 & 0 \\
0 & 0 & 0 & 1
 \end{smallmatrix}  \right)$; hence $(\pi(r)\pi(s))^2 \neq (\pi(s)\pi(r))^2$. Let 

\begin{align*}
R_1 &= \left\{\left(\begin{smallmatrix} 1 & * &* & *  \\
0 & 1 & * & * \\
0 & 0 &2 & * \\
0 & 0 & 0 & 2\end{smallmatrix}\right) \in \F_3^{4\times 4} \right\},& R &= \pi(R_1) \cap \oc \ni \pi(r); \\
S_1 &= \left\{\left(\begin{smallmatrix} 
1 & * &* & *  \\
0 & 2 & * & * \\
0 & 0 &2 & * \\
0 & 0 & 0 & 1\end{smallmatrix}\right) \in \F_3^{4\times 4} \right\},& S &= \pi(S_1) \cap \oc \ni \pi(s).
\end{align*}
Then $R \coprod S \hookrightarrow \oc$ is a decomposable subrack,  
and $\oc$ is of type D.
\epf

%\newpage

\end{document}